\documentclass[11pt]{article}
\usepackage{amssymb,amsmath,color} % font used for R in Real numbers
 \allowdisplaybreaks

\catcode`!=11
\let\!int\int \def\int{\displaystyle\!int}
\let\!lim\lim \def\lim{\displaystyle\!lim}
\let\!sum\sum \def\sum{\displaystyle\!sum}
\let\!sup\sup \def\sup{\displaystyle\!sup}
\let\!inf\inf \def\inf{\displaystyle\!inf}
\let\!cap\cap \def\cap{\displaystyle\!cap}
\let\!max\max \def\max{\displaystyle\!max}
\let\!min\min \def\min{\displaystyle\!min}
\let\!frac\frac \def\frac{\displaystyle\!frac}
\catcode`!=12

\let\oldsection\section
\renewcommand\section{\setcounter{equation}{0}\oldsection}

\allowdisplaybreaks

\newtheorem{thm}{Theorem}[section]
\newtheorem{pro}{Proposition}[section]
\newtheorem{lem}{Lemma}[section]
\newtheorem{re}{Remark}[section]
\newtheorem{cor}{Corollary}[section]

\begin{document}
%%%%%%%%%%%%%%%%%%%% title %%%%%%%%%%%%%%%%%%%%%%%%%%%%%%%%%%%%%%%%%%%%%%%%%%
\title{LOCAL WELL-POSEDNESS OF  PRANDTL EQUATIONS FOR COMPRESSIBLE FLOW IN TWO SPACE VARIABLES}
\author{{\bf Ya-Guang Wang,\quad Feng
Xie}\\[1mm]
\small Department of Mathematics, and LSC-MOE,\\[1mm]
\small Shanghai Jiao Tong University,
Shanghai 200240, P.R.China\\[1mm]
{\bf Tong Yang}\\[1mm]
\small Department of Mathematics,
City University of Hong Kong,\\[1mm]
\small Tat Chee Avenue, Kowloon, Hong Kong
}
\date{}
\maketitle

\noindent{\bf Abstract:} In this paper, we consider the local
well-posedness of the Prandtl boundary layer equations
that describe the behavior of boundary layer in the small viscosity limit of the compressible isentropic Navier-Stokes equations with non-slip boundary condition. Under the strictly monotonic assumption on
 the tangential velocity in the normal variable, we apply
 the Nash-Moser-H\"{o}rmander
iteration scheme and further develop the energy method introduced in
\cite{AWXY} to obtain the well-posedness of the equations locally in time. 

%In contrast to the work \cite{AWXY} for the incompressible Prandtl equations, one key point of this work is that one needs %to construct the first approximate solution as the starting point of iteration satisfying the compressible Prandtl equations %with faster decay errors. As a consequence, we deduce the well-posedness of the incompressible Prandtl equations in two %space variables with non-uniform outer flows.

\footnotetext[2]{\indent  {\it E-mail address:} ygwang@sjtu.edu.cn(Y.G. Wang),
tzxief@sjtu.edu.cn (F. Xie), matyang@cityu.edu.hk(T. Yang)}

\vskip 2mm


\vskip 2mm


\section{Introduction}
In this paper, we study the well-posedness of the compressible Prandtl boundary layer equations that are derived in the small viscosity limit from the compressible
isentropic Naiver-Stokes equations with non-slip boundary condition. 
Note that the Prandtl equations describe the behavior of the  characteristic boundary layer in the leading order.
% (Refer to the Appendix for derivations in details). 
Denote by $\mathbb{T}\times\mathbb{R}^+=\{(x, \eta)|x\in\mathbb{R}/\mathbb{Z}, 0\leq \eta<+\infty\}$ the periodic spatial domain, 
and let $u(t,x,\eta)$ and $v(t,x,\eta)$ be the tangential and normal velocity 
components in the boundary layer. Consider the following compressible Prandtl equations with $(x,\eta)\in \mathbb{T}\times\mathbb{R}^+$,
\begin{align} \label{1.1}
\left\{
\begin{array}{ll}
u_t+uu_x+vu_{\eta}-\frac{1}{\bar{\rho}(t,x)}\partial_{\eta}^2u+P_x=0,\\
\partial_x(\bar{\rho}u)+\partial_{\eta}(\bar{\rho}v)=-\bar{\rho}_t,
\end{array}
\right.
\end{align}
with the initial data
\begin{align}
\label{1.1I} u(t,x,\eta)|_{t=0}=u_0(x,\eta),
\end{align}
and  the boundary and the far-field conditions
\begin{align}
\label{1.2} u(t,x,\eta)|_{\eta=0}=0,\quad v(t,x,\eta)|_{\eta=0}=0,\quad
\lim_{\eta\rightarrow+\infty}u(t,x,\eta)=U(t,x).
\end{align}
Here, $\bar{\rho}(t,x)$ and $U(t,x)$ are the traces on the  boundary $\{y=0\}$  of  the density  and the tangential velocity of the outer Euler flow that
satisfy the Bernoulli's law
\begin{equation}\label{ber}
U_t+UU_x+P_x=0,
\end{equation}
with $P(t,x)$ being the trace of the enthalpy of the outer Euler flow.

It is well-known that the leading order characteristic boundary layer for the
incompressible Navier-Stokes equations 
with non-slip boundary condition is described by the classical Prandtl
equations that were proposed by Prandtl  \cite{P} in 1904.
%
%From that time on it attracts much attention from mathematicians.
%However, the well-posedness (or ill-posedness) of solutions
%to the Prandtl equations are far from being understood due to the difficulties in mathematical analysis
%\cite{AWXY,CS,E,EE,GD,GN,GR,GN1,O,OS,XY}.
Under the monotone
assumption on the tangential velocity in the normal direction, Oleinik firstly
obtained the local existence of classical solutions  in the two spatial
dimension
by using the Crocco transformation, cf.  \cite{O}. This result
together with some other extensions in this direction are presented in
Oleinik-Samokhin's classical
book \cite{OS}. Recently, this well-posedness result was re-proved by using an energy method  in the framework of
 Sobolev spaces in \cite{AWXY} and \cite{MW1} independently. On
the other hand, by imposing an additional  
favorable condition on the pressure, 
a global in time weak solution was obtained in
\cite{XZ}.

When the monotonicity condition is violated, seperation of the boundary
layer is well expected and observed. For this, E-Engquist constructed a finite
time blowup solution to the Prandtl
equations in \cite{EE}. Recently, when the background shear flow has a non-degenerate critical point, some interesting ill-posedness (or instability) phenomena of solutions to both
the linear and nonlinear Prandtl equations around
the shear flow are studied, cf.
\cite{GD,GN,GR,GN1}. All these results show that the monotone assumption on the tangential velocity is very important for well-posedness 
except in the framework of analytic functions
studied in \cite{CS} and some other references with generalization.

\iffalse

By applying the multi-scale analysis for the solutions to the compressible isentropic Navier-Stokes equations with characteristic boundary conditions, the behavior of boundary
layers can be described by  the equations (\ref{1.1}). It is noted that the
characteristic boundary conditions for compressible Navier-Stokes
equations also comprise the no-slip boundary condition and the Navier-slip
boundary conditions with friction coefficients depending on the
small viscosity. One can refer to \cite{MW, WXW} for the physical description and
derivations in details. To complete the presentation, we give the derivation of (\ref{1.1}) in details in the appendix.

\fi

This paper aims to obtain the local well-posedness of the problem (\ref{1.1})-(\ref{1.2}) for the
compressible Prandtl equations  in some weighted
Sobolev spaces. To state the main results,  we first give 
the following
assumptions on the initial data.\\

\underline{\bf Main assumptions (H) on the initial data}:
\begin{enumerate}
\item[(H1)] For a fixed integer $k_0\geq 9$, the initial data $u_0(x,\eta)$ satisfies the compatibility
condition of the problem  (\ref{1.1})-(\ref{1.2}) up to order $4k_0+2$;
\item[(H2)] Monotone condition $\partial_{\eta}u_0(x,\eta)\geq \frac{\sigma_0}{(1+\eta)^{\gamma+2}}>0$ holds for all $x\in\mathbb{T}$ and $\eta\geq 0$ with some 
 positive constant $\sigma_0$ and a positive integer $\gamma\geq 2$;
\item[(H3)] $\|(1+\eta)^{\gamma+\alpha_2}D^{\alpha}(u_0(x,\eta)-U(0,x))\|_{L^2(\mathbb{T}\times\mathbb{R}^+)}\leq C_0$, where $D^{\alpha}=\partial_x^{\alpha_1}\partial_{\eta}^{\alpha_2}$ with $\alpha=(\alpha_1, \alpha_2)$ and $|\alpha|=\alpha_1+\alpha_2\leq 4k_0+2$;
\item[(H4)] $\|(1+\eta)^{\gamma+2+\alpha_2}D^{\alpha}\partial_{\eta}u_{0}\|_{L^{\infty}(\mathbb{T}\times\mathbb{R}^+)}\leq \frac{1}{\sigma_0}$, for $|\alpha|\leq 3k_0$.
\end{enumerate}

Denote by $V(t,x)$ the trace of $\partial_yu_2^E$ on $\{y=0\}$ for the normal velocity $u_2^E$ of Euler outer flow.  From the conservation of mass in the Euler equations, we have 
$$\partial_t \bar{\rho}(t,x)+\partial_x(\bar{\rho}(t,x)U(t,x))+\bar{\rho}(t,x)V(t,x)=0.$$
Here, we have used the fact that $u_2^E(t,x,y)|_{y=0}=0$. Thus, from the problem (\ref{1.1})-(\ref{1.2}), we know that the normal velocity $v(t,x,\eta)$ can be represented by
\begin{align}
\label{1.1A} v(t,x,\eta)=V(t,x)\eta+\frac{1}{\bar{\rho}(t,x)}\int_0^{\eta}\partial_x(\bar{\rho}(t,x)(U(t,x)-u(t,x, \tilde{\eta})))d\tilde{\eta}.
\end{align}

The main result of this paper can be stated  as follows.
 \begin{thm}
 \label{MAIN}
Suppose that the outer Euler flow is smooth for $0\le t\le T_0$,
the density $\bar{\rho}(t,x)$ has both positive lower and upper bounds, and the
Sobolev norm $H^{s}([0,T_0]\times \mathbb{R})$ of $(\bar{\rho}, U, V)$ is bounded for a suitably large integer 
$s$, moreover, the Main Assumption (H) on the initial data
$u_0(x,\eta)$ is satisfied. Then there exists  $0<T\le T_0$, such
that the initial boundary value problem (\ref{1.1})-(\ref{1.2}) has
a unique classical solutions $(u, v)$ satisfying
\begin{align}
\label{1.5}
\sum_{|m_1|+[(m_2+1)/2]\leq k_0}
\|\langle\eta\rangle^l\partial_{(t,x)}^{m_1}\partial_{\eta}^{m_2}(u-U)\|_{L^2([0,T]\times\mathbb{T}\times\mathbb{R}^+)}<+\infty,
\end{align}
for a fixed $l>\frac 12$ depending only on $\gamma$ given in (H)  with $\langle\eta\rangle=(1+\eta)$, and
\begin{align}
\label{1.5-1}
\sum_{|m_1|+[(m_2+1)/2]\leq k_0-1} \sup\limits_{\eta\in \mathbb{R}^+}
\|\partial_{(t,x)}^{m_1}\partial_{\eta}^{m_2}(v-V\eta)(\cdot, \eta)\|_{L^2([0,T]\times\mathbb{T})}<+\infty.
\end{align}
 \end{thm}

\begin{re} (1)
When the outer Euler flow density $\bar{\rho}(t,x)$ is a positive constant,
  the system (\ref{1.1}) is reduced to the classical incompressible Prandtl equations. Thus the analysis  in this paper works also for the classical incompressible Prandtl equations with general far-field condition and initial data satisfying the Main Assumption (H). Note that
 the case with a uniform outer flow with slightly different
assumption on the initial data was studied in \cite{AWXY}.

(2) It is straightforward to verify that the set of the initial data satisfying the Main Assumption (H) is not empty
because it contains the functions
with polynomial decay in  $\eta$. 
%A similar assumption of polynomial decreasing is also required in \cite{MW1}.
\end{re}

Now, let us give some comments on the analysis in this paper. 
In principle, we will apply the approach of \cite{AWXY} to study the problem
(\ref{1.1})-(\ref{1.2}). There are several crucial differences
between the system (\ref{1.1}) and classical incompressible Prandtl
equations. Firstly, the normal velocity $v$ contains the linearly
increasing part $V(t,x)\eta$ in $\eta$, consequently, in estimating the
solution to the linearized problem, we  need to study the
 conormal estimates.
Secondly, the divergence free condition in the classical Prandtl
system is now replaced by an inhomogeneous equation in (\ref{1.1}).
Moreover,  the far-field state is not uniform 
so that  the shear flow is no longer an
exact solution to the compressible Prandtl equations (\ref{1.1}).
Therefore, to apply the Nash-Moser-H\"{o}mander iteration
scheme used in \cite{AWXY} for the nonlinear problem
(\ref{1.1})-(\ref{1.2}), we need to construct a suitable zero-th
order approximate solution with suitable
error estimate. And the construction is given in subsection 4.1
in three steps.

%Due to the facts again that the far-field state is not uniform and the divergence equation is inhomogenous, the
%techniques used in \cite{AWXY} for the energy estimates should be revised correspondingly. In addition, the domain considered in
%\cite{AWXY} is half plane, we consider our problem in the periodic
%domain $\mathbb{T}\times\mathbb{R}^+=\{(x, \eta)|x\in\mathbb{R}/\mathbb{Z}, 0\leq \eta<+\infty\}$.

Finally, the rest of the paper is organized as follows. We will
 first introduce some
weighted Sobolev spaces and give some preliminaries in Section 2. The
well-posedness  of the linearized compressible Prandtl
equations is given in Section 3. In Section 4, we introduce the
Nash-Moser-H\"omander iteration scheme, and construct the first
approximate solution as the starting point of iteration. Then
  the local existence and uniqueness of solution to the
nonlinear problem of the compressible Prandtl equations are proved. 

% Finally, we derive the boundary layer equations (\ref{1.1})-(\ref{1.2}) for the compressible viscous flows in Appendix.

 \section{Preliminaries}
In this section, we will introduce some weighted Sobolev spaces and
norms for later use. To simplify the notations, we denote by
$\partial_{\tau}^m$ the summation of tangential derivatives
$\partial_{\tau}^m=\partial_{t}^{m_1}\partial_{x}^{m_2}$ for all
$m=(m_1,m_2)\in\mathbb{N}^2, |m|=m_1+m_2$. Denote
$\partial_{\tau}^{\alpha}$ by $Z_1^{\alpha}$,
$(\eta\partial_{\eta})^{\alpha}$ by $Z^{\alpha}_2$  and
$Z^{m}=Z_1^{m_1}Z^{m_2}_2$ for further simplification. Define
\begin{align*}
\|f\|_{H_{co,l}^{k_1,k_2}}=\left(\sum_{0\leq m_1\leq k_1, 0\leq
m_2\leq k_2}
\|\langle\eta\rangle^lZ^{m_1}_1Z^{m_2}_2f\|_{L^2([0,T]\times\mathbb{T}\times
\mathbb{R}^+)}\right)^{1/2},
\end{align*}
and
\begin{align*}
\|f\|_{D_{co,l}^{k_1,k_2}}=\left(\sum_{0\leq m_1\leq k_1, 0\leq
m_2\leq k_2}\sup_{\eta\ge 0}
\|\langle\eta\rangle^lZ^{m_1}_1Z^{m_2}_2f(\cdot,
\eta)\|_{L^2([0,T]\times\mathbb{T})}\right)^{1/2},
\end{align*}
\begin{align*}
\|f\|_{C_{co,l}^{k_1,k_2}}=\left(\sum_{0\leq m_1\leq k_1, 0\leq
m_2\leq k_2}\sup_{(t,x)\in [0,T]\times \mathbb{T}}
\|\langle\eta\rangle^lZ^{m_1}_1Z^{m_2}_2f(t,x,\cdot)\|_{L^{2}(\mathbb{R}^+)}\right)^{1/2}.
\end{align*}
Denote 
$$\|f\|_{H^k_{co,l}}=\sum_{k_1+k_2=k}\|f\|_{H^{k_1,k_2}_{co,l}}.$$
The function spaces $D^k_{co,l}$ and $C^k_{co,l}$ can be defined
similarly.  Since the conormal operator
$Z^{m}$ does not communicate with the normal
derivative operator $\partial_{\eta}$, the following estimate of
commutator is frequently used,
\begin{align}
\label{CE} \|[Z^{m}, \partial_{\eta}]f\|_{L^2_l}\lesssim
\|\partial_{\eta}f\|_{H^{m-1}_{co,l}}.
\end{align}
Here and after $0<a \lesssim  b$
%$a \lesssim(\gtrsim) b$
means that there exists a
uniform constant $C>0$ such that $a \leq Cb$.\\
The following weighted Sobolev spaces are also used frequently.
Denote by $\partial_{\eta}^k$ the $k$-th normal derivative, for any
given $k_1,k_2\in\mathbb{N}, l\in\mathbb{R}^+$ and $0<T<+\infty$, 
set
\begin{align*}
\|f\|_{\mathcal{B}^{k_1,k_2}_{l}}=\left(\sum_{0\leq |m|\leq k_1,
0\leq n\leq k_2}\|\langle\eta\rangle^lZ^m\partial^n_{\eta}f\|^2_{L^2([0,T]
\times\mathbb{T}\times\mathbb{R}^+)}\right)^{1/2},
\end{align*}
\begin{align*}
\|f\|_{\tilde{\mathcal{B}}^{k_1,k_2}_{l}}=\left(\sum_{0\leq
|m|\leq k_1, 0\leq n\leq k_2}\|\langle\eta\rangle^lZ^m
\partial^n_{\eta}f\|^2_{L^{\infty}([0,T];L^2(\mathbb{T}\times\mathbb{R}^+))}\right)^{1/2},
\end{align*}
\begin{align*}
\|f\|_{\mathcal{A}_l^m}=\left(\sum_{|m_1|+[(m_2+1)/2]\leq|m|}
\|\langle\eta\rangle^{l}Z^{m_1}\partial_{\eta}^{m_2}f\|^2_{L^2([0,T]
\times\mathbb{T}\times\mathbb{R}^+)}\right)^{1/2}.
\end{align*}
It is straightforward to verify that
\begin{align*}
\mathcal{A}^m_{l}=\bigcap_{m_1+[(m_2+1)/2]\leq m}\mathcal{B}_{l}^{m_1,m_2}.
\end{align*}
We also define
\begin{align*}
\|f\|_{\mathcal{D}^{m}_l}=\left(\sum_{k_1+[(k_2+1)/2]\leq m}
\|\langle\eta\rangle^lZ^{k_1}\partial^{k_2}_{\eta}f\|^2_{L^{\infty}_{\eta}(L^2_{t,x})}\right)^{1/2},
\end{align*}
and
\begin{align*}
\|f\|_{\mathcal{C}^{m}_l}=\left(\sum_{k_1+[(k_2+1)/2]\leq m}
\|\langle\eta\rangle^lZ^{k_1}\partial^{k_2}_{\eta}f\|^2_{L^{\infty}_{t,x}(L^2_{\eta})}\right)^{1/2}.
\end{align*}
In addition, the homogeneous norms
$\|\cdot\|_{\dot{\mathcal{A}}_l^m},
\|\cdot\|_{\dot{\mathcal{C}}_l^m},
\|\cdot\|_{\dot{\mathcal{D}}_l^m}$ correspond to the summation over
$1\leq |m_1|+[(m_2+1)/2]\leq |m|$.\\
For $1\leq p\leq +\infty$, we will use
$\|f\|_{L_l^p(\mathbb{T}\times\mathbb{R}^+)}=
\|\langle\eta\rangle^lf\|_{L^p(\mathbb{T}\times\mathbb{R}^+)}$. It
is direct to show the following Sobolev type embeddings,
\begin{align}
\label{2.1} \|f\|_{\mathcal{C}^m_l}\leq
C_s\|f\|_{\mathcal{A}_l^{m+2}},\quad \|f\|_{\mathcal{D}^m_l}\leq
C_s\|f\|_{\mathcal{A}_{l+1}^{m+1}}.
\end{align}
Moreover, for any $l\geq 0$ and $m\geq 2$, the space
$\mathcal{A}_l^m$ is continuously embedded into
$\mathcal{C}_b^{m-2}$ which is the space of $(m-2)-$th order continuously
differentiable functions with bounded derivatives. And the
following Morse-type inequalities hold.
\begin{lem}
\label{l2.1} For any proper functions $f$ and $g$, we have
\begin{align*}
\|fg\|_{\mathcal{A}_l^m}\leq
C_m\left\{\|f\|_{\mathcal{A}_l^m}\|g\|_{L^{\infty}}+
\|f\|_{L^{\infty}}\|g\|_{\dot{\mathcal{A}}_l^m}\right\},
\end{align*}
and
\begin{align*}
\|fg\|_{\mathcal{A}_l^m}\leq
C_m\left\{\|f\|_{\mathcal{C}_l^m}\|g\|_{\mathcal{D}^0_0}+
\|f\|_{\mathcal{C}_l^0}\|g\|_{\dot{\mathcal{D}}_0^m}\right\}.
\end{align*}
Similar inequalities hold in the norms
$\|\cdot\|_{H_{co, l}^{m}}, \|\cdot\|_{\mathcal{B}_l^{m_1,m_2}}, \|\cdot\|_{\mathcal{C}_l^m}$
and $\|\cdot\|_{\mathcal{D}_l^m}$. Here, $C_m>0$ is a constant
depending only on $m$.
\end{lem}

These results can be obtained  similarly  as those given in \cite{ME}.

\section{Well-posedness of linearized system}
The  strategy to prove the main result, Theorem \ref{MAIN}, is
to apply an
 iteration scheme to construct a sequence of approximate
solution sequences, and then to show these approximate solutions
converge in some suitable weighted Sobolev space. 
Since there is a loss of
regularity, the Nash-Moser-H\"omander iteration scheme is used
for this purpose. In
this section, we study the well-posedness of the linearized equations
and obtain the required
energy estimates of solutions to the
linearized equations for the Nash-Moser-H\"omander iteration. \\
Let $(\tilde{u},
\tilde{v})$ be a smooth background state satisfying the following conditions.
$$
\partial_{\eta}\tilde{u}(t,x,\eta)>0,\quad
\partial_{x}(\bar{\rho}\tilde{u})+\partial_{\eta}(\bar{\rho}\tilde{v})
=-\bar{\rho}_t,\quad \tilde{u}|_{\eta=0}=\tilde{v}|_{\eta=0}=0,
\quad \lim_{\eta\rightarrow +\infty}\tilde{u}=U(t,x).
$$
Here $\tilde{v}$ is given by
$$
\tilde{v}=V(t,x)\eta+\frac{1}{\bar{\rho}(t,x)}\int_0^{\eta}
\partial_x(\bar{\rho}(t,x)(U(t,x)-\tilde{u}))d\tilde{\eta}\triangleq V(t,x)\eta+\bar{v}.
$$
It extracts the linear increasing part $V(t,x)\eta$ by introducing
the new function $\bar{v}$. The linearized problem of
(\ref{1.1})-(\ref{1.2}) around $(\tilde{u}, \tilde{v})$ can be
written as 
\begin{align}
\label{3.1} \left\{
\begin{array}{ll}
u_t+\tilde{u}u_x+\tilde{v}\partial_{\eta}u+u\tilde{u}_{x}
+\tilde{u}_{\eta}v-\frac{1}{\bar{\rho}}u_{\eta\eta}=f,\\
\partial_{\eta}(\bar{\rho}v)+\partial_x(\bar{\rho}u)=0,\\
u|_{\eta=0}=v|_{\eta=0}=0,\quad\lim_{\eta\rightarrow+\infty}u(t,x,\eta)=0,\\
u|_{t=0}=0.
\end{array}
\right.
\end{align}
Similar to \cite{AWXY}, by introducing the  transformation
$$
\omega(t,x,\eta)=\left(\frac{\bar{\rho}u}{\partial_{\eta}\tilde{u}}\right)_{\eta}(t,x,\eta),
$$
then for classical solutions, from (\ref{3.1}) we know that $w$
satisfies the following problem in $\{t>0, x\in \mathbb{T},
\eta>0\}$:
\begin{align}
\label{3.2} \begin{cases}
\begin{array}{ll}
\omega_t+(\tilde{u}\omega)_x+(\eta V \omega)_{\eta}+(\bar{v}\omega)_{\eta}-\frac{2}{\bar{\rho}}\left(\omega\frac{\partial_{\eta}^2\tilde{u}}{\partial_{\eta}\tilde{u}}\right)_{\eta}
+\left(\xi\int_0^{\eta}\omega(t,x,\tilde{\eta})d\tilde{\eta}\right)_{\eta}-\frac{\bar{\rho}_t}{\bar{\rho}}\omega-\frac{1}{\bar{\rho}}
\omega_{\eta\eta}=\tilde{f}_{\eta},\\
\frac{1}{\bar{\rho}}\left(\omega_{\eta}+2\omega\frac{\partial_{\eta}^2\tilde{u}}{\partial_{\eta}\tilde{u}}\right)|_{\eta=0}=-\tilde{f}|_{\eta=0},\\
\omega|_{t=0}=0,
\end{array}
\end{cases}
\end{align}
where
$$
\xi=\frac{(\partial_t+\tilde{u}\partial_x+\tilde{v}\partial_{\eta}-\frac{1}{\bar{\rho}}\partial_{\eta}^2)\tilde{u}_{\eta}}{\tilde{u}_{\eta}}
-\frac{\tilde{u}\bar{\rho}_x}{\bar{\rho}}\triangleq\xi_1-\frac{\tilde{u}\bar{\rho}_x}{\bar{\rho}},\quad\quad
\tilde{f}=\frac{\bar{\rho}f}{\tilde{u}_{\eta}}.
$$
To simplify the presentation, we use the notations:
\begin{align}\label{lambda}
\lambda_{k_1,k_2}=\|\tilde{u}-U\|_{\mathcal{B}^{k_1,k_2}_{l}}+
\|Z^{k_1}\partial_{\eta}^{k_2}\bar{v}\|_{L^{\infty}_{\eta}(L^2_{t,x})}
+\|Z^{k_1}\partial_{\eta}^{k_2}\chi\|_{L^{\infty}_{\eta}(L^{2}_{t,x})}+\|\xi_1\|_{\mathcal{B}^{k_1,k_2}_{l}},
\end{align}
with
$
\chi=\frac{\partial^2_{\eta}\tilde{u}}{\partial_{\eta}\tilde{u}}.
$
Set
\begin{align}\label{lambda-k}
\lambda_k=\sum_{k_1+[(k_2+1)/2]\leq k}\lambda_{k_1,k_2}.
\end{align}
Similar to \cite{AWXY}, we have the following energy estimates of the solution to the
problem (\ref{3.2}).
\begin{thm}
\label{Thm3.1} Suppose that the outer Euler flow
$(\bar{\rho}(t,x), U(t,x), V(t,x))\in H^s(\mathbb{R}^2_+)$, for $s$ suitably large, and $\bar{\rho}(t,x)$ has
uniform lower positive bound. Moreover,  for a given positive $k$, the compatibility condition
for the problem (\ref{3.2}) holds up to order $k$. Then for any fixed
$l>1/2$, we have 
\begin{align}
\label{3.3} \|\omega\|_{\mathcal{A}^k_l}\leq
C_1(\lambda_4)\|\tilde{f}\|_{\mathcal{A}^k_l}+C_2(\lambda_4)\lambda_k\|\tilde{f}\|_{\mathcal{A}^3_l},
\end{align}
with $C_1(\cdot)$ and $C_2(\cdot)$ being two smooth functions in their arguments.
\end{thm}

As we mentioned in the introduction, the main difference of the linear problem \eqref{3.2} from the one studied in \cite{AWXY} is that there is a linear growth term $\eta V$ in the equation of \eqref{3.2}.
Hence, we can  not obtain the estimates of tangential derivatives directly as in \cite{AWXY}. Similar to \cite{MW}, we will study the linearized problem \eqref{3.2} in  some conormal space.  First, we have
\begin{lem}
($L^2$-estimate) Under the assumptions in Theorem 3.1, there exists a positive constant $C$ such that
\begin{align}\label{lemma3.1}
\frac{d}{dt}\|\omega\|^2_{L^2_l(\mathbb{T}\times\mathbb{R}^+)}+\|\omega_{\eta}\|^2_{L^2_l(\mathbb{T}\times\mathbb{R}^+)}
\leq C(\lambda_{3,1}+1)\|\omega\|^2_{L^2_l(\mathbb{T}\times\mathbb{R}^+)}+C\|\tilde{f}\|^2_{L^2_l(\mathbb{T}\times\mathbb{R}^+)}.
\end{align}
\end{lem}

\begin{pf}
Multiplying (\ref{3.2}) by $\langle\eta\rangle^{2l}\omega$ and
integrating it over $\mathbb{T}\times\mathbb{R}^+$, we obtain
\begin{align}
\label{Basic} 
\frac{d}{dt}\|\omega\|^2_{L^2_l(\mathbb{T}\times\mathbb{R}^+)} & +2\int_{\mathbb{T}\times\mathbb{R}^+}\langle\eta\rangle^{2l}\omega\left\{(\tilde{u}\omega)_x+(\eta V \omega)_{\eta}+(\bar{v}\omega)_{\eta}-\frac{2}{\bar{\rho}}\left(\omega\frac{\partial_{\eta}^2\tilde{u}}{\partial_{\eta}\tilde{u}}\right)_{\eta}\right.\nonumber\\
&\left.+\left(\xi\int_0^{\eta}\omega(t,x,\tilde{\eta})d\tilde{\eta}\right)_{\eta}-\frac{\bar{\rho}_t}{\bar{\rho}}\omega-\frac{1}{\bar{\rho}}\omega_{\eta\eta}-\tilde{f}_{\eta}\right\}dxd\eta=0.
\end{align}
It is straightforward to obtain
\begin{align*}
\int_{\mathbb{T}\times\mathbb{R}^+}\{(\tilde{u}\omega)_x+(\bar{v}\omega)_{\eta}\}\langle\eta\rangle^{2l}\omega dxd\eta
& =\int_{\mathbb{T}\times\mathbb{R}^+}\langle\eta\rangle^{2l}\frac{\omega^2}2
(\frac{-\bar{\rho}_t-\bar{\rho}V-\tilde{u}\bar{\rho}_x}{\bar{\rho}})-l\bar{v}\omega^2\langle\eta\rangle^{2l-1}dxd\eta\\
\lesssim&\left(1+\|\tilde{u}\|_{L^{\infty}(\mathbb{T}\times\mathbb{R}^+)}+\|\bar{v}\|_{L^{\infty}(\mathbb{T}\times\mathbb{R}^+)}
\right)\|\omega\|^2_{L^2_l(\mathbb{T}\times\mathbb{R}^+)},
\end{align*}
and
\begin{align*}
&\int_{\mathbb{T}\times\mathbb{R}^+}(\eta V\omega)_{\eta}\langle\eta\rangle^{2l}\omega dxd\eta
=\int_{\mathbb{T}\times\mathbb{R}^+}V\omega^2(\frac12\langle\eta\rangle^{2l}-l\eta\langle\eta\rangle^{2l-1})dxd\eta
\lesssim \|\omega\|^2_{L^2_l(\mathbb{T}\times\mathbb{R}^+)}.
\end{align*}

On the other hand, by integration by parts and using the boundary condition given in
(\ref{3.2}) we get 
\begin{align*}
\int_{\mathbb{T}\times\mathbb{R}^+}\left\{-\frac2{\bar{\rho}}(\omega\chi)_{\eta}-\frac{1}{\bar{\rho}}\omega_{\eta\eta}-\tilde{f}_{\eta}\right\}
\langle\eta\rangle^{2l}\omega dxd\eta
=\int_{\mathbb{T}\times\mathbb{R}^+}\left\{\frac2{\bar{\rho}}(\omega\chi)+\frac{1}{\bar{\rho}}\omega_{\eta}+\tilde{f}\right\}
(\langle\eta\rangle^{2l}\omega)_{\eta}dxd\eta,
\end{align*}
where the right hand side can be estimated as follows.
\begin{align*}
&\int_{\mathbb{T}\times\mathbb{R}^+}\frac{2}{\bar{\rho}}\omega\chi (\langle\eta\rangle^{2l}\omega)_{\eta}dxd\eta
\lesssim \|\chi\|_{L^{\infty}(\mathbb{T}\times\mathbb{R}^+)}(\|\omega\|_{L_l^2(\mathbb{T}\times\mathbb{R}^+)}^2
+\|\omega\|_{L_l^2(\mathbb{T}\times\mathbb{R}^+)}\|\partial_{\eta}\omega\|_{L_l^2(\mathbb{T}\times\mathbb{R}^+)}),
\end{align*}
\begin{align*}
\int_{\mathbb{T}\times\mathbb{R}^+}\frac{1}{\bar{\rho}}\omega_{\eta}(2l\langle\eta\rangle^{2l-1}\omega+\langle\eta\rangle^{2l}\omega_{\eta})
\gtrsim \|\partial_{\eta}\omega\|^2_{L^2_l(
\mathbb{T}\times\mathbb{R}^+)}-\|\omega\|^2_{L^2_l(\mathbb{T}\times\mathbb{R}^+)},
\end{align*}
and
\begin{align*}
\int_{\mathbb{T}\times\mathbb{R}^+}\tilde{f}(2l\langle\eta\rangle^{2l-1}\omega+\langle\eta\rangle^{2l}\omega_{\eta})
\lesssim\|\tilde{f}\|_{L^2_l(\mathbb{T}\times\mathbb{R}^+)}(\|\omega\|_{L^2_l(\mathbb{T}\times\mathbb{R}^+)}+\|\partial_{\eta}\omega\|_{L^2_l(\mathbb{T}\times\mathbb{R}^+)}).
\end{align*}
%And \begin{align*}
%\left|\int_{\mathbb{T}\times\mathbb{R}^+}(-\frac{\bar{\rho}_t{\bar{\rho}}\omega)\langle\eta\rangle^{2l}\omega dxd\eta\right|\lesssim\|\omega\|%^2_{L^2_l(\mathbb{T}\times\mathbb{R}^+)}.
%\end{align*}
Denote by
\begin{align*}
\int_{\mathbb{T}\times\mathbb{R}^+}\left((\xi_1-\frac{\tilde{u
}\bar{\rho}_x}{\bar{\rho}})\int_{0}^{\eta}\omega(t,x,\tilde{\eta})d\tilde{\eta}\right)_{\eta}\langle\eta\rangle^{2l}\omega dxd\eta\triangleq H^1+H^2.
\end{align*}
As $l>1/2$, by integration by parts, it follows
\begin{align*}
|H^1|\lesssim\|\xi_1\|_{L_x^{\infty}(L^2_{\eta,l})}(\|\omega\|^2_{L^2_l(\mathbb{T}\times\mathbb{R}^+)}
+\|\omega\|_{L^2_l(\mathbb{T}\times\mathbb{R}^+)}\|\partial_{\eta}\omega\|_{L^2_l(\mathbb{T}\times\mathbb{R}^+)}),
\end{align*}
and
\begin{align*}
|H^2|\leq&|\int_{\mathbb{T}\times\mathbb{R}^+}\frac{\bar{\rho}_x(\tilde{u}-U)}{\bar{\rho}}\int_0^{\eta}\omega(t,x,\tilde{\eta})d\tilde{\eta}
(\langle\eta\rangle^{2l}\omega)_{\eta}dxd\eta|\\
&+|\int_{\mathbb{T}\times\mathbb{R}^+}\left(\frac{\bar{\rho}_xU}{\bar{\rho}}\int_{0}^{\eta}\omega(t,x,\tilde{\eta})d\tilde{\eta}\right)_{\eta}
\langle\eta\rangle^{2}\omega dxd\eta|\\
\lesssim&\|\omega\|^2_{L^2_l(\mathbb{T}\times\mathbb{R}^+)}
+\|\tilde{u}-U\|_{L^{\infty}_x(L^2_{\eta,l})}(\|\omega\|^2_{L^2_l(\mathbb{T}\times\mathbb{R}^+)}
+\|\omega\|_{L^2_l(\mathbb{T}\times\mathbb{R}^+)}\|\partial_{\eta}\omega\|_{L^2_l(\mathbb{T}\times\mathbb{R}^+)}).
\end{align*}

\iffalse

It should be remarked that in the estimates above we use the facts that $\bar{\rho}$ has positive lower and upper bounds and $|U, V, \bar{\rho}_t, \bar{\rho}_x|$ are bounded without any emphasis.
By the classical Sobolev inequalities, we have
$$
\|\tilde{u}, \bar{v},
\chi\|_{L^{\infty}(\mathbb{T}\times\mathbb{R}^+)},
\qquad\|\xi_1,(\tilde{u}-U)\|_{L^{\infty}_x(L^2_{\eta, l})}\lesssim
(1+\lambda_{3,1}).
$$
For simplicity, here and after, we always assume
the $L^{\infty}$ and $L^2_{x,t}$ norms of the known functions $\bar{\rho}, U, V$ and
its derivatives of all order are bounded. Moreover, $\bar{\rho}(t,x)$ has a
positive lower bound due to the assumptions in Theorem \ref{MAIN}.
It is noted that these assumptions
are reasonable, at least in local time.

\fi

 Thus, from \eqref{Basic} we obtain 
\begin{align*}
\frac{d}{dt}\|\omega\|^2_{L^2_l(\mathbb{T}\times\mathbb{R}^+)}+\|\omega_{\eta}\|^2_{L^2_l(\mathbb{T}\times\mathbb{R}^+)}
\leq C(1+\lambda_{3,1})\|\omega\|^2_{L^2_l(\mathbb{T}\times\mathbb{R}^+)}+C\|\tilde{f}\|^2_{L^2_l(\mathbb{T}\times\mathbb{R}^+)},
\end{align*}
by noting $$
\|\tilde{u}, \bar{v},
\chi\|_{L^{\infty}(\mathbb{T}\times\mathbb{R}^+)}+\|\xi_1\|_{L^{\infty}_x(L^2_{\eta, l})}+
\|\tilde{u}-U\|_{L^{\infty}_x(L^2_{\eta, l})}\lesssim
(1+\lambda_{3,1}).
$$
It completes the proof of the estimate \eqref{lemma3.1}.

\end{pf}
\begin{lem}
(Estimates of conormal derivatives) Under the assumptions in Theorem 3.1, for any fixed
$T>0$, there exists a positive constant $C>0$ such that
\begin{align}\label{3.6}
 &\frac{d}{dt}\|\omega\|^2_{H^{m}_{co,l}}+\|\omega_{\eta}\|^2_{H^{m}_{co,l}}\nonumber\\
 \leq&
C(1+\lambda^2_{3,1})\|\omega\|^2_{H^{m}_{co,l}}+C(\|\tilde{f}\|^2_{H^{m}_{co,l}}+\|\tilde{f}_{\eta}\|^2_{H^{m-1}_{co,l}}+(\lambda^2_{m-1,1}+\lambda^2_{m,0}+1)\|\omega\|^2_{\mathcal{B}^{3,1}_{l}})
\end{align}
holds for $t\in[0,T]$.
\end{lem}

\begin{pf}
The proof is divided into four steps.

(1) Applying the conormal derivative operator $Z^{m}$ on the
equation in (\ref{3.2}), multiplying the resulting equation by
$\langle\eta\rangle^{2l}Z^{m}\omega$ and
integrating it over $\mathbb{T}\times\mathbb{R}^+$,  it follows
\begin{align}
\label{3.7A} 
\frac{d}{dt}\|Z^m\omega\|^2_{L^2_l(\mathbb{T}\times\mathbb{R}^+)} & +2\int_{\mathbb{T}\times\mathbb{R}^+}\langle\eta\rangle^{2l}Z^m \omega Z^m\left\{(\tilde{u}\omega)_x+(\eta V \omega)_{\eta}+(\bar{v}\omega)_{\eta}-\frac{2}{\bar{\rho}}\left(\omega\frac{\partial_{\eta}^2\tilde{u}}{\partial_{\eta}\tilde{u}}\right)_{\eta}\right.\nonumber\\
&\left.+\left(\xi\int_0^{\eta}\omega(t,x,\tilde{\eta})d\tilde{\eta}\right)_{\eta}-\frac{\bar{\rho}_t}{\bar{\rho}}\omega-\frac{1}{\bar{\rho}}\omega_{\eta\eta}-\tilde{f}_{\eta}\right\}dxd\eta=0.
\end{align}

Now, let us estimate each term of \eqref{3.7A}.
Denote
$$\begin{cases}
I_1=\int_{\mathbb{T}\times\mathbb{R}^+}Z^m[(\tilde{u}\omega)_x+(\bar{v}\omega)_{\eta}]\langle\eta\rangle^{2l}Z^{m}\omega dxd\eta,\\
I_2=-\int_{\mathbb{T}\times\mathbb{R}^+}\langle\eta\rangle^{2l}Z^m \omega Z^m\left(\frac{2}{\bar{\rho}}(\omega\frac{\partial_{\eta}^2\tilde{u}}{\partial_{\eta}\tilde{u}})_\eta+\frac{1}{\bar{\rho}}\omega_{\eta\eta}+\tilde{f}_{\eta}\right)dxd\eta,\\
I_3=\int_{\mathbb{T}\times\mathbb{R}^+}\langle\eta\rangle^{2l}Z^m \omega Z^m\left(\xi\int_0^{\eta}\omega(t,x,\tilde{\eta})d\tilde{\eta}\right)_{\eta}dxd\eta.
\end{cases}
$$

(2) { \it Estimate of $I_1$.} Obviously, we have
\begin{align*}
I_1=&\int_{\mathbb{T}\times\mathbb{R}^+}Z^m[(\tilde{u}\omega)_x+(\bar{v}\omega)_{\eta}]\langle\eta\rangle^{2l}Z^{m}\omega dxd\eta\\
=&\int_{\mathbb{T}\times\mathbb{R}^+}Z^m[(-\frac{\bar{\rho}_t+\bar{\rho}V}{\bar{\rho}})\omega+(-\frac{\bar{\rho}_x\tilde{u}}{\bar{\rho}})\omega+\tilde{u}\omega_x
+\bar{v}\omega_{\eta}]\langle\eta\rangle^{2l}Z^m\omega dxd\eta\\
\triangleq& I_1^1+I_1^2+I_1^3+I_1^4,
\end{align*}
and 
\begin{align*}
|I_1^1|=\left|\sum_{m_1+m_2=m}C_m^{m_1}\int_{\mathbb{T}\times\mathbb{R}^+}\langle\eta\rangle^{2l}\left(Z^{m_1}(-\frac{\bar{\rho}_t+\bar{\rho}V}{\bar{\rho}})\right)(Z^{m_2}\omega)(Z^m\omega) dxd\tau\right|
\lesssim\|\omega\|^2_{H^m_{co,l}}, 
\end{align*}
\begin{align*}
|I_1^2|
=&\left|\sum_{m_1+m_2=m}C_m^{m_1}\int_{\mathbb{T}\times\mathbb{R}^+}Z^{m_1}\left[(\frac{\bar{\rho}_x(\tilde{u}-U(x,t))}
{\bar{\rho}})+\frac{\bar{\rho}_xU}{\bar{\rho}}\right]Z^{m_2}\omega
\langle\eta\rangle^{2l}Z^m\omega dxd\eta\right|\\
\lesssim&\|\tilde{u}-U\|_{H^m_{co,l}}\|\omega\|_{L^{\infty}}\|\omega\|_{H^m_{co,l}}+(1+\|\tilde{u}\|_{L^{\infty}})\|\omega\|^2_{H^m_{co,l}}.
\end{align*}

On the other hand, one has
\begin{align*}
I_1^3=&\sum_{m_1+m_2=m,m_2<m}C_m^{m_1}\int_{\mathbb{T}\times\mathbb{R}^+}\langle\eta\rangle^{2l}
(Z^{m_1}\tilde{u})
(Z^{m_2}\omega_x)
(Z^m\omega) dxd\eta\\
&+\int_{\mathbb{T}\times\mathbb{R}^+}\langle\eta\rangle^{2l}\tilde{u}(Z^{m}\omega_x)(Z^m\omega) dxd\eta
\triangleq I_1^{3a}+I_1^{3b},
\end{align*}
and
\begin{align*}
I_1^4=&\sum_{m_1+m_2=m,m_2<m}C_m^{m_1}\int_{\mathbb{T}\times\mathbb{R}^+}
\langle\eta\rangle^{2l}
(Z^{m_1}\bar{v})(Z^{m_2}\omega_{\eta})(Z^m\omega) dxd\eta\\
&+\int_{\mathbb{T}\times\mathbb{R}^+}\langle\eta\rangle^{2l}\bar{v}(Z^{m}\omega_{\eta})(Z^m\omega) dxd\eta
\triangleq I_1^{4a}+I_1^{4b}.
\end{align*}

Note that $I_1^{3a}$ can be estimated similarly as  $I_1^2$, and 
\begin{align*}
|I_1^{4a}|\lesssim\|Z^m\omega\|_{L^2_l}(\|Z^m\bar{v}\|_{L^{\infty}_{\eta}(L^2_{x})}\|\omega_{\eta}\|_{L^2_{\eta}(L^{\infty}_{x})}
+\|Z\bar{v}\|_{L^{\infty}}\|Z^{m-1}\omega_{\eta}\|_{L^2_l}).
\end{align*}
By using $\partial_{\eta}(\bar{\rho}\bar{v})+\partial_{x}(\bar{\rho}\tilde{u})=-\bar{\rho}_t-\bar{\rho}V$, integration by parts and the commutator estimate (\ref{CE}), we obtain
\begin{align*}
|I_1^{3b}+I_1^{4b}|\lesssim (1+\|\tilde{u}\|_{L^{\infty}}+\|\bar{v}\|_{L^{\infty}})
\|Z^m\omega\|^2_{L^2_l(\mathbb{T}\times\mathbb{R}^+)}+\|\bar{v}\|_{L^{\infty}}\|\partial_{\eta}\omega\|_{H^{m-1}_{co,l}}\|Z^m\omega\|_{L^2_l(\mathbb{T}\times\mathbb{R}^+)}.
\end{align*}
Moreover, the definition of the operator $Z_2$ gives that
\begin{align*}
\int_{\mathbb{T}\times\mathbb{R}^+}Z^m(\eta V\omega)_{\eta}\langle\eta\rangle^{2l}Z^m\omega dxd\eta
=\int_{\mathbb{T}\times\mathbb{R}^+}Z^m(V\omega+V\eta \omega_{\eta})\langle\eta\rangle^{2l}Z^m\omega dxd\eta
\lesssim \|\omega\|^2_{H^m_{co,l}}.
\end{align*}

(3) { \it Estimate of $I_2$.} First, by using the boundary condition given in \eqref{3.2}, we have
\begin{align*}
I_2=&\int_{\mathbb{T}\times\mathbb{R}^+}(\langle\eta\rangle^{2l}Z^m\omega)_{\eta}
Z^m\left(\frac{2}{\bar{\rho}}\omega\frac{\partial_{\eta}^2\tilde{u}}{\partial_{\eta}\tilde{u}}+\frac{1}{\bar{\rho}}\omega_{\eta}+\tilde{f}\right)dxd\eta\\
&-
\int_{\mathbb{T}\times\mathbb{R}^+}\langle\eta\rangle^{2l}Z^m\omega
[Z^m, \partial_\eta]
\left(\frac{2}{\bar{\rho}}\omega\frac{\partial_{\eta}^2\tilde{u}}{\partial_{\eta}\tilde{u}}+\frac{1}{\bar{\rho}}\omega_{\eta}+\tilde{f}\right)dxd\eta\\
\triangleq& I_2^1+I_2^2+I_2^3+I_2^4,
\end{align*}
with $I_2^4$ being the terms involving $\tilde f$,
$$
I_2^1=\int_{\mathbb{T}\times\mathbb{R}^+}(\langle\eta\rangle^{2l}Z^m\omega)_{\eta}
Z^m\left(\frac{2}{\bar{\rho}}\omega\frac{\partial_{\eta}^2\tilde{u}}{\partial_{\eta}\tilde{u}}\right)dxd\eta
-
\int_{\mathbb{T}\times\mathbb{R}^+}\langle\eta\rangle^{2l}Z^m\omega
[Z^m, \partial_\eta]
\left(\frac{2}{\bar{\rho}}\omega\frac{\partial_{\eta}^2\tilde{u}}{\partial_{\eta}\tilde{u}}\right)dxd\eta,$$
$$I_2^2=\int_{\mathbb{T}\times\mathbb{R}^+}(\langle\eta\rangle^{2l}Z^m\omega)_{\eta}
Z^m\left(\frac{1}{\bar{\rho}}\omega_{\eta}\right)dxd\eta,$$
and
$$I_2^3=-
\int_{\mathbb{T}\times\mathbb{R}^+}\langle\eta\rangle^{2l}Z^m\omega
[Z^m, \partial_\eta]
\left(\frac{1}{\bar{\rho}}\omega_{\eta}\right)dxd\eta.$$

It is straightforward to show that
\begin{align*}
|I^1_2|\leq&\left|\int_{\mathbb{T}\times\mathbb{R}^+}Z^m(\frac{2}{\bar{\rho}}\chi
\omega)\langle\eta\rangle^{2l}(Z^m\omega_{\eta}+\frac{2l}{\langle\eta\rangle}Z^m\omega+[Z^m,\partial_{\eta}]\omega)dxd\eta\right|
+\|Z^m\omega\|_{L^2_l}\|\frac2{\bar{\rho}}(\chi \omega)_{\eta}\|_{H^{m-1}_{co,l}}\\
\lesssim&\|\frac{2}{\bar{\rho}}\chi \omega\|_{H^m_{co,l}}(\|Z^m\omega_{\eta}\|_{L^2_l}+\|Z^m\omega\|_{L^2_l}+\|\partial_{\eta}\omega\|_{H^{m-1}_{co,l}})
+\|Z^m\omega\|_{L^2_l}\|\frac2{\bar{\rho}}(\chi \omega)_{\eta}\|_{H^{m-1}_{co,l}}\\
\lesssim&(\|\chi\|_{\l^{\infty}}\|\omega\|_{H^{m}_{co,l}}+\|\chi\|_{D^{m}_{co,0}}\|\omega\|_{L^{2}_{\eta,l}(L^{\infty}_x)})
(\|\omega\|_{H^{m}_{co,l}}+\|\partial_{\eta}\omega\|_{H^{m}_{co,l}})\\
+\|Z^m&\omega\|_{L^2_l}(\|\partial_{\eta}\chi\|_{D^{m-1,0}_{co,0}}\|\omega\|_{L^{2}_{\eta,l}(L^{\infty}_x)}
+\|\partial_{\eta}\chi\|_{L^{\infty}}\|\omega\|_{H^{m-1}_{co,l}}+\|\chi\|_{D^{m-1,0}_{co,0}}\|\omega_{\eta}\|_{L^{2}_{\eta,l}(L^{\infty}_x)}
+\|\chi\|_{L^{\infty}}\|\omega_{\eta}\|_{H^{m-1}_{co,l}}),
\end{align*}
\begin{align*}
I_2^2=&\int_{\mathbb{T}\times\mathbb{R}^+}
(Z^{m}\frac{\omega_{\eta}}{\bar{\rho}})(2l\langle\eta\rangle^{2l-1}Z^m\omega+\langle\eta\rangle^{2l}(Z^m\omega)_{\eta})dxd\eta\\
\gtrsim&\frac{1}{2}\|Z^m\omega_{\eta}\|^2_{L^2_l}-\|\omega\|^2_{H^m_{co,l}}-\|\omega_{\eta}\|^2_{H^{m-1}_{co,l}},
\end{align*}
and
\begin{align*}
|I_2^3|
\leq\|\partial_{\eta}(\frac{\omega_{\eta}}{\bar{\rho}})\|_{H^{m-1}_{co,l}}\|Z^m\omega\|_{L^2_l}.
\end{align*}
From the equation (\ref{3.2}), we have
\begin{align*}
&\|\partial_{\eta}(\frac{\omega_{\eta}}{\bar{\rho}})\|_{H^{m-1}_{co,l}}\\
=&\|\omega_t+(\tilde{u}\omega)_x+(\eta V \omega)_{\eta}+(\bar{v}\omega)_{\eta}-\frac{2}{\bar{\rho}}\left(\omega\frac{\partial_{\eta}^2\tilde{u}}{\partial_{\eta}\tilde{u}}\right)_{\eta}
+\left(\xi\int_0^{\eta}\omega(t,x,\tilde{\eta})d\tilde{\eta}\right)_{\eta}-\frac{\bar{\rho}_t}{\bar{\rho}}\omega-\tilde{f}_{\eta}\|_{H^{m-1}_{co,l}}.
\end{align*}
The terms on the right hand side of the above equation
can be estimated as 
\begin{align*}
\|\omega_t+(\eta V \omega)_{\eta}-\frac{\bar{\rho}_t}{\rho}\omega\|_{H^{m-1}_{co,l}}\lesssim \|\omega\|_{H^{m}_{co,l}},
\end{align*}
\begin{align*}
&\|(\tilde{u}\omega)_x+(\bar{v}\omega)_{\eta}-\frac{2}{\bar{\rho}}\left(\omega\frac{\partial_{\eta}^2\tilde{u}}{\partial_{\eta}\tilde{u}}\right)_{\eta}\|_{H^{m-1}_{co,l}}\\
&=\|-\left(\frac{\bar{\rho}_t+\bar{\rho}V+\bar{\rho}_x\tilde{u}}{\bar{\rho}}
+\frac{2}{\bar{\rho}}\left(\frac{\partial_{\eta}^2\tilde{u}}{\partial_{\eta}\tilde{u}}\right)_{\eta}\right)
\omega
+\tilde{u}\omega_x+\left(\bar{v}-
\frac{2}{\bar{\rho}}\frac{\partial_{\eta}^2\tilde{u}}{\partial_{\eta}\tilde{u}}
\right)\omega_{\eta}\|_{H^{m-1}_{co,l}}\\
& \leq (1+\|\tilde{u}\|_{L^{\infty}}+\|\chi_\eta\|_{L^\infty})\|\omega\|_{H^{m}_{co,l}}
+(\|\bar{v}\|_{L^{\infty}}+\|\chi\|_{L^\infty})\|\omega_{\eta}\|_{H^{m-1}_{co,l}}\\
& +(\|\tilde{u}-U\|_{H^{m-1}_{co,l}}+\|\chi_\eta\|_{H^{m-1}_{co,l}})
(\|\omega\|_{L^{\infty}}+\|\omega_x\|_{L^{\infty}})
+(\|\bar{v}\|_{D^{m,0}_{co,l}}+\|\chi\|_{D^{m,0}_{co,l}})\|\omega_{\eta}\|_{L^2_{\eta,l}(L^{\infty}_x)},
\end{align*}
and
\begin{align*}
&\|\left(\xi\int_0^{\eta}\omega(t,x,\tilde{\eta})d\tilde{\eta}\right)_{\eta}\|_{H^{m-1}_{co,l}}\\
&\le (\|\xi_1\|_{H^{m-1}_{co,l}}+\|{\tilde u}-U\|_{H^{m-1}_{co,l}})
\|\omega\|_{L^{\infty}}+\|\xi\|_{L^{\infty}}\|\omega\|_{H^{m-1}_{co,l}}\\
&+\|(\xi_1)_{\eta}, ({\tilde u}-U)_\eta\|_{H^{m-1}_{co,l}}\|\omega\|_{L^2_{\eta,l}(L^{\infty}_x)}+\|(\xi_1)_{\eta},  ({\tilde u}-U)_\eta\|_{L^2_{\eta,l}(L^{\infty}_x)}\|\omega\|_{H^{m-1}_{co,l}}.
\end{align*}

(4) { \it Estimate of $I_3$.} 
Decompose $I_3$ into
$$I_3=I_3^1+I_3^2,$$
with
$$I_3^1=
\int_{\mathbb{T}\times\mathbb{R}^+}\langle\eta\rangle^{2l}Z^m \omega Z^m\left(\xi_1\int_0^{\eta}\omega(t,x,\tilde{\eta})d\tilde{\eta}\right)_{\eta}dxd\eta,
$$
and
$$I_3^2=-
\int_{\mathbb{T}\times\mathbb{R}^+}\langle\eta\rangle^{2l}Z^m \omega Z^m\left(\frac{{\tilde u}{\bar \rho}_x}{\bar \rho}\int_0^{\eta}\omega(t,x,\tilde{\eta})d\tilde{\eta}\right)_{\eta}dxd\eta.
$$
For $l>1/2$, we have
\begin{align*}
|I^1_3|\leq&\left|\int_{\mathbb{T}\times\mathbb{R}^+}
Z^m\left(\xi_1\int_0^{\eta}\omega d\tilde{\eta}\right)(\langle\eta\rangle^{2l}Z^m\omega)_{\eta}dxd\eta\right|+
\left|\int_{\mathbb{T}\times\mathbb{R}^+}
[Z^m,\partial_{\eta}]\left(\xi_1\int_0^{\eta}\omega d\tilde{\eta}\right)(\langle\eta\rangle^{2l}Z^m\omega)dxd\eta\right|\\
\lesssim&(\|\xi_1\|_{L^2_{\eta,l}(L^{\infty}_x)}\|\omega\|_{H^{m}_{co,l}}+\|Z^m\xi_1\|_{L^2_l}
\|\omega\|_{L^2_{\eta,l}(L^{\infty}_x)})(\|Z^m\omega\|_{L^2_l}+\|\omega_{\eta}\|_{H^m_{co,l}})\\
&+\|Z^m\omega\|_{L^2_{l}}
\|\partial_{\eta}(\xi_1\int_0^{\eta}\omega d\tilde{\eta})\|_{H^{m-1}_{co,l}}\\
\lesssim&(\|\xi_1\|_{L^2_{\eta,l}(L^{\infty}_x)}\|\omega\|_{H^{m}_{co,l}}+\|Z^m\xi_1\|_{L^2_l}
\|\omega\|_{L^2_{\eta,l}(L^{\infty}_x)})(\|Z^m\omega\|_{L^2_l}+\|\omega_{\eta}\|_{H^m_{co,l}})\\
&+\|Z^m\omega\|_{L^2_l}(\|\xi_1\|_{H^{m-1}_{co,l}}\|\omega\|_{L^{\infty}_l}+\|\xi_1\|_{L^{\infty}_l}\|\omega\|_{H^{m-1}_{co,l}}\\
&+\|(\xi_1)_{\eta}\|_{H^{m-1}_{co,l}}\|\omega\|_{L^2_{\eta,l}(L^{\infty}_x)}+\|(\xi_1)_{\eta}\|_{L^2_{\eta,l}(L^{\infty}_x)}\|\omega\|_{H^{m-1}_{co,l}}),
\end{align*}
and
\begin{align*}
|I_3^2|=&\left|\int_{\mathbb{T}\times\mathbb{R}^+}Z^m\left(\frac{\bar{\rho}_x(\tilde{u}-U+U)}{\bar{\rho}}\int_0^{\eta}
\omega(t,x,\tilde{\eta})d\tilde{\eta}\right)_{\eta}
\langle\eta\rangle^{2l}Z^m\omega dxd\eta\right|\\
\lesssim&(\|Z^m(\tilde{u}-U)\|_{L^2_l}\|\omega\|_{L^2_{\eta,l}(L^{\infty}_x)}+\|\tilde{u}-U\|_{L^2_{\eta,l}(L^{\infty}_x)}\|\omega\|_{H^m_{co,l}})
(\|Z^m\omega\|_{L^2_l}+\|\omega_{\eta}\|_{H^m_{co,l}})\\
&+\|Z^m\omega\|_{L^2_l}(\|\partial_{\eta}(\tilde{u}-U)\|_{L^{\infty}_l}\|\omega\|_{H^{m-1}_{co,l}}
+\|\partial_{\eta}(\tilde{u}-U)\|_{H^{m-1}_{co,l}}\|\omega\|_{L^2_{\eta,l}(L^{\infty}_x)}\\
&+\|\tilde{u}-U\|_{H^{m-1}_{co,l}}\|\omega\|_{L^{\infty}}+\|\tilde{u}-U\|_{L^{\infty}}\|\omega\|_{H^{m-1}_{co,l}}+\|Z^m\omega\|_{L^2_l}).
\end{align*}

Summarizing the above estimates, it follows
\begin{align}
\label{3.6}
&\frac{d}{dt}\|\omega\|^2_{H^{m}_{co,l}}+\|\omega_{\eta}\|^2_{H^{m}_{co,l}}\nonumber\\
\lesssim& \|\tilde{f}\|^2_{H^{m}_{co,l}}+\|\tilde{f}_{\eta}\|^2_{H^{m-1}_{co,l}}+(1+\lambda^2_{3,1})\|\omega\|^2_{H^{m}_{co,l}}
+(\lambda^2_{m,0}+\lambda^2_{m-1,1}+1)\|\omega\|^2_{\mathcal{B}^{3,1}_{l}},
\end{align}
where we have used the inequalities
\begin{align*}
\|\omega, \omega_x\|_{L^{\infty}_l(\mathbb{T}\times\mathbb{R}^+)},\quad
\|\omega, \omega_{\eta}\|_{L^2_{\eta,l}(L^{\infty}_x)}\leq
C\|\omega\|_{\mathcal{B}^{3,1}_{l}}.
\end{align*}
And this completes the proof of the lemma.
\end{pf}

\begin{re} Similar to the above proof, one can obtain
\begin{align}
\label{3.7} \|\omega\|^2_{\mathcal{B}^{m,1}_{l}}\leq
C(\|\tilde{f}\|^2_{H^m_{co,l}}+\|\tilde{f}_{\eta}\|^2_{H^{m-1}_{co,l}}),\quad 0\leq m\leq 3.
\end{align}
When $m=0$,  the term $\|\tilde{f}_{\eta}\|^2_{H^{m-1}_{co,l}}$ 
is not in (\ref{3.7}). By combining (\ref{3.6}), (\ref{3.7}) and using Gronwall's inequality, we get
\begin{align}
\label{3.6A} \|\omega\|^2_{\mathcal{B}^{m,1}_{l}}\leq
C(\|\tilde{f}\|^2_{\mathcal{B}^{m,0}_{l}}+\|\tilde{f}\|^2_{\mathcal{B}^{m-1,1}_{l}})+(1+\lambda^2_{m-1,1}+\lambda^2_{m,0})
(\|\tilde{f}\|^2_{\mathcal{B}^{3,0}_{l}}+\|\tilde{f}\|^2_{\mathcal{B}^{2,1}_{l}}).
\end{align}
\end{re}

We are now ready to prove Theorem \ref{Thm3.1}.\\

{\it Proof of Theorem \ref{Thm3.1}:}
It remains to  estimate the higher order normal derivatives.

From the equation given in (\ref{3.2}), we have
\begin{align}
\label{3.8}
\frac{1}{\bar{\rho}}\omega_{\eta\eta}=\omega_t+(\tilde{u}\omega)_x+(\eta V \omega)_{\eta}+(\bar{v}\omega)_{\eta}-\frac{2}{\bar{\rho}}\left(\omega\frac{\partial_{\eta}^2\tilde{u}}{\partial_{\eta}\tilde{u}}\right)_{\eta}
+\left(\xi\int_0^{\eta}\omega(t,x,\tilde{\eta})d\tilde{\eta}\right)_{\eta}-\frac{\bar{\rho}_t}{\bar{\rho}}\omega-\tilde{f}_{\eta}.
\end{align}
Applying the conormal operator $Z^{m}$ to the above equation (\ref{3.8})
gives
\begin{equation}
\label{3.12}
\begin{array}{ll}
\|\omega\|_{\mathcal{B}^{m,2}_{l}}\leq&
\|\bar{\rho}\omega_t\|_{\mathcal{B}^{m,0}_{l}}+\|\bar{\rho}((\tilde{u}\omega)_x+(\bar{v}\omega)_{\eta})\|_{\mathcal{B}^{m,0}_{l}}
+\|\bar{\rho}(\eta V\omega)_{\eta}\|_{\mathcal{B}^{m,0}_{l}}+\|2\chi
\omega\|_{\mathcal{B}^{m,1}_{l}}\\
&+\|\bar{\rho}\partial_{\eta}(\xi\int_{0}^{\eta}\omega(t,x,\tilde{\eta})d\tilde{\eta})\|_{\mathcal{B}^{m,0}_{l}}
+\|\bar{\rho}_t\omega\|_{\mathcal{B}^{m,0}_{l}}
+\|\bar{\rho}\partial_{\eta}\tilde{f}\|_{\mathcal{B}^{m,0}_{l}}.
\end{array}
\end{equation}
We  estimate each term on the right  hand side of the above inequality. By using Lemma
\ref{l2.1}, we get
\begin{align*}
\|\bar{\rho}\omega_t\|_{\mathcal{B}^{m,0}_{l}}+\|(\eta V\omega)_{\eta}\|_{\mathcal{B}^{m,0}_{l}}+\|\bar{\rho}_t\omega\|_{\mathcal{B}^{m,0}_{l}}\lesssim \|\omega\|_{\mathcal{B}^{m+1,0}_{l}}.
\end{align*}
Obviously, it holds
\begin{align*}
\|\bar{\rho}((\tilde{u}\omega)_x+(\bar{v}\omega)_{\eta})\|_{\mathcal{B}^{m,0}_{l}}
=\|(\bar{\rho}_t+\bar{\rho}V+\bar{\rho}_x\tilde{u})\omega-\bar{\rho}\tilde{u}\omega_x-\bar{\rho}\bar{v}\omega_{\eta}\|_{\mathcal{B}^{m,0}_{l}},
\end{align*}
where
\begin{align*}
\|(\bar{\rho}_t+\bar{\rho}V)\omega\|_{\mathcal{B}^{m,0}_{l}}\lesssim\|\omega\|_{\mathcal{B}^{m,0}_{l}},
\end{align*}
\begin{align*}
\|\bar{\rho}\bar{v}\omega_{\eta}\|_{\mathcal{B}^{m,0}_{l}}\lesssim \|\bar{v}\|_{L^{\infty}}\|\omega_{\eta}\|_{\mathcal{B}^{m,0}_{l}}+\|\bar{v}\|_{\mathcal{D}^{m,0}_{co,0}}\|\omega_{\eta}\|_{L^2_{\eta,l}(L^{\infty}_x)},
\end{align*}
and
\begin{align*}
\|\bar{\rho}\tilde{u}\omega_x\|_{\mathcal{B}^{m,0}_{l}}=&\|\bar{\rho}(\tilde{u}-U)\omega_x+\bar{\rho}U\omega_x\|_{\mathcal{B}^{m,0}_{l}}\\
\lesssim&\|\tilde{u}-U\|_{\mathcal{B}^{m,0}_{l}}\|\omega_x\|_{L^{\infty}}+(\|\tilde{u}-U\|_{L^\infty}+1)\|\omega\|_{\mathcal{B}^{m+1,0}_{l}}.
\end{align*}
The term $\|\bar{\rho}_x\tilde{u}\omega\|_{\mathcal{B}^{m,0}_{l}}$ can be estimated similarly.
Moreover, we have
\begin{align*}
\|2\chi \omega\|_{\mathcal{B}_l^{m,1}}\lesssim \|\chi\|_{L^{\infty}}\|\omega\|_{\mathcal{B}_l^{m,1}}+\|\chi\|_{\mathcal{D}^{m,1}_{0}}\|\omega\|_{L^{2}_{\eta,l}(L^{\infty}_x)}.
\end{align*}
And 
\begin{align*}
\|\bar{\rho}\partial_{\eta}(\xi\int_{0}^{\eta}\omega(t,x,\tilde{\eta})d\tilde{\eta})\|_{\mathcal{B}^{m,0}_{l}}
=\|\bar{\rho}(\xi_{\eta}\int_{0}^{\eta}\omega(t,x,\tilde{\eta})d\tilde{\eta}+\xi \omega)\|_{\mathcal{B}^{m,0}_{l}},
\end{align*}
where
\begin{align*}
\|\bar{\rho}(\xi_{\eta}\int_{0}^{\eta}\omega(t,x,\tilde{\eta})d\tilde{\eta}\|_{\mathcal{B}^{m,0}_{l}}
\lesssim&\|\xi_{1\eta}\|_{\mathcal{B}^{m,0}_{l}}\|\omega\|_{L^2_{\eta,l}(L^{\infty}_x)}+\|\xi_{1\eta}\|_{L^{2}_{\eta,l}(L^{\infty}_{x,t})}\|\omega\|_{\mathcal{B}^{m,0}_{l}}\\
&+\|\tilde{u}-U\|_{\mathcal{B}_l^{m,1}}\|\omega\|_{L^2_{\eta,l}(L^{\infty}_x)}+\|\tilde{u}\|_{L^{\infty}_l}\|\omega\|_{\mathcal{B}_l^{m,1}},
\end{align*}
and
\begin{align*}
\|\xi \omega\|_{\mathcal{B}^{m,0}_{l}}\lesssim
\|\xi_1\|_{\mathcal{B}^{m,0}_{l}}\|\omega\|_{L^{\infty}_{\eta,l}}+\|\xi_1\|_{L^{\infty}_l}\|\omega\|_{\mathcal{B}^{m,0}_{l}}
+\|\tilde{u}-U\|_{\mathcal{B}^{m,0}_{l}}\|\omega\|_{L^{\infty}}+(\|\tilde{u}\|_{L^{\infty}}+1)\|\omega\|_{\mathcal{B}^{m,1}_{l}}.
\end{align*}
Plugging the above estimates into \eqref{3.12} yields
\begin{align*}
\|\omega\|_{\mathcal{B}^{m,2}_{l}}\lesssim\lambda_{3,1}(\|\omega\|_{\mathcal{B}^{m,1}_{l}}+\|\omega\|_{\mathcal{B}^{m+1,0}_{l}})
+\lambda_{m,1}\|\omega\|_{\mathcal{B}^{3,1}_{l}}+\|\tilde{f}\|_{\mathcal{B}^{m,1}_{l}}.
\end{align*}

Next, for any fixed $n\ge 3$, applying the differential operator $Z^m\partial^{n-2}_{\eta}$ on the equation (\ref{3.8}) gives
\begin{align*}
&\|\omega\|_{\mathcal{B}^{m,n}_{l}}\\
=&\|\bar{\rho}(\omega_t+(\tilde{u}\omega)_x+(\eta V \omega)_{\eta}+(\bar{v}\omega)_{\eta}-\frac{2}{\bar{\rho}}\left(\omega\frac{\partial_{\eta}^2\tilde{u}}{\partial_{\eta}\tilde{u}}\right)_{\eta}
+\left(\xi\int_0^{\eta}\omega(t,x,\tilde{\eta})d\tilde{\eta}\right)_{\eta}-\frac{\bar{\rho}_t}{\bar{\rho}}\omega-\tilde{f}_{\eta})\|_{\mathcal{B}^{m,n-2}_{l}},
\end{align*}
where
\begin{align*}
\|\bar{\rho}(\omega_t+(\eta V \omega)_{\eta}-\frac{\bar{\rho}_t}{\bar{\rho}}\omega\|_{\mathcal{B}^{m,n-2}_{l}}\lesssim \|\omega\|_{\mathcal{B}^{m+1,n-2}_{l}},
\end{align*}
\begin{align*}
\|\bar{\rho}(\tilde{u}\omega)_x+\bar{\rho}(\bar{v}\omega)_{\eta}\|_{\mathcal{B}^{m,n-2}_{l}}\lesssim& \|\omega\|_{\mathcal{B}^{m,n-2}_{l}}+\|\tilde{u}-U\|_{\mathcal{B}^{m,n-2}_{l}}\|\omega_x\|_{L^{\infty}}
+(\|\tilde{u}-U\|_{L^{\infty}}+1)\|\omega\|_{\mathcal{B}^{m+1,n-2}_{l}}\\
&+\|\bar{v}\|_{L^{\infty}}\|\omega_{\eta}\|_{\mathcal{B}^{m,n-2}_{l}}+\|\bar{v}\|_{\mathcal{D}^{m,n-2}_0}\|\omega_{\eta}\|_{L^2_{\eta,l}(L^{\infty}_x)},
\end{align*}
\begin{align*}
\|2(\omega\chi)_{\eta}\|_{\mathcal{B}^{m,n-2}_{l}}\lesssim \|\omega\|_{\mathcal{B}^{m,n-1}_{l}}\|\chi\|_{L^{\infty}}+\|\chi\|_{\mathcal{D}^{m,n-1}_0}\|\omega\|_{L^{2}_{\eta,l}(L^{\infty}_x)},
\end{align*}
and
\begin{align*}
&\|\left(\xi\int_0^{\eta}\omega(t,x,\tilde{\eta})d\tilde{\eta}\right)_{\eta}\|_{\mathcal{B}^{m,n-2}_{l}}
=\|\xi_{\eta}\int_0^{\eta}\omega(t,x,\tilde{\eta})d\tilde{\eta}+\xi \omega\|_{\mathcal{B}^{m,n-2}_{l}}\\
\lesssim&\|\xi_{\eta}\|_{\mathcal{B}^{m,n-2}_{l}}\|\omega\|_{L^2_{\eta,l}(L^{\infty}_x)}+\|\xi_{\eta}\|_{L^{\infty}_l}\|\omega\|_{\mathcal{B}^{m,n-2}_{l}}
+\|\xi\|_{\mathcal{B}^{m,n-2}_{l}}\|\omega\|_{L^{\infty}_{\eta,l}}+\|\xi\|_{L^{\infty}_l}\|\omega\|_{\mathcal{B}^{m,n-2}_{l}}.
\end{align*}
Thus, we obtain
\begin{align*}
\|\omega\|_{\mathcal{B}^{m,n}_{l}}\lesssim (\lambda_{3,1}+1)(\|\omega\|_{\mathcal{B}^{m+1,n-2}_{l}}+\|\omega\|_{\mathcal{B}^{m,n-1}_{l}})+\lambda_{m,n-1}\|\omega\|_{\mathcal{B}^{3,1}_{l}}+\|\tilde{f}\|_{\mathcal{B}^{m,n-1}_{l}}.
\end{align*}
By induction on $n$, we conclude the estimate (\ref{3.3}).
And this completes the proof of Theorem 3.1.

\section{Iteration scheme and convergence}
Based on the energy estimate \eqref{3.3} on the solution to the linearized
equations obtained in the previous
 section, we now study the well-posedness of the nonlinear problem (\ref{1.1}) by using a suitable linear
iteration scheme. From \eqref{3.3},  there is a loss of regularity in the
solutions to the linearized problem \eqref{3.1} with respect to both of the background states and initial data.
Hence, as in \cite{AWXY}, we apply the Nash-Moser-H\"omander iteration
scheme.   As we explained in Section 1, 
we do not have the divergence free condition, 
and the far-field state is not uniform.
Thus,  the shear flow is no longer the special
exact solution to the compressible Prandtl equations (\ref{1.1}) in contrast to the incompressible problem studied in \cite{AWXY}. Thus, to start the Nash-Moser-H\"{o}mander iteration, we need to construct a proper zero-th
order approximate solution  satisfying the nonlinear compressible Prandtl equations with enough decay in $\eta$. The construction
will be given Subsection 4.1. Then, in Subsection 4.2 we 
present the Nash-Moser-H\"omander iteration
scheme for the problem (\ref{1.1}).  The estimates of the approximate solutions are obtained in Subsection 4.3. In Subsection 4.4, we conclude the convergence of iteration for the existence and uniqueness of the solution to the nonlinear problem (\ref{1.1})-(\ref{1.2}).

\subsection{The Zero-th order approximate solution}
In this subsection, we construct the initial approximate solution to 
in the following three subsections.

\subsubsection{Compatibility conditions and initial data}

Set
$$
u=U(t,x)+\bar{u}, \quad v=V(t,x)\eta+\bar{v}.
$$
By using the Bernoulli law \eqref{ber}, from (\ref{1.1}) we know that $(\bar{u}, \bar{v})$ satisfies
\begin{align}
\label{APP01}
\left\{
\begin{array}{ll}
\bar{u}_t+U_x\bar{u}+U\bar{u}_x+\bar{u}\bar{u}_x+(\bar{v}+V\eta)\bar{u}_{\eta}-\frac{1}{\bar{\rho}}\partial_{\eta}^2\bar{u}=0,\\
\partial_{\eta}(\bar{\rho}\bar{v})+\partial_x(\bar{\rho}\bar{u})=0,\\
\bar{u}(t,x,\eta)|_{t=0}=u_0-U(0,x), \quad \bar{v}|_{\eta=0}=0.
\end{array}
\right.
\end{align}
Denote
$$
\bar{u}^j(x,\eta)=\partial_t^j\bar{u}(t,x,\eta)|_{t=0}, \quad
\bar{v}^j(x,\eta)=\partial_t^j\bar{v}(t,x,\eta)|_{t=0}.
$$
From the compatibility condition of (\ref{APP01}),
$\{\bar{u}^j, \bar{v}^j\}_{j\leq k_0}$ is in turn given explicitly by
$u_0(x,\eta),  U(0,z)$ and $V(0,x)$.\\
We  define the first approximate solution
$(\bar{u}, \bar{v})$ of (\ref{APP01}) as follows.
\begin{align}
\label{4.12}
u^a(t,x,\eta)=\sum_{j=0}^{k_0}\frac{t^j}{j!}\bar{u}^j(x,\eta),
\quad v^a(t,x,\eta)=-\frac{1}{\bar{\rho}}\int_0^{\eta}(\bar{\rho}u^a)_x(t,x,\tilde{\eta})d\tilde{\eta}.
\end{align}
From (H3) and (H4) in the Main Assumptions (H), it follows that
\begin{align}
\label{APP06}
\max_{0\le t\le T}\|\langle\eta\rangle^{\gamma+\alpha_2}D^{\alpha}u^a(t,\cdot)\|_{L^2({\mathbb{T}\times\mathbb{R}^+})}\leq CC_0, \quad |\alpha|\leq 2k_0,
\end{align}
for a fixed $T>0$, where $C$ depends on $\sigma_0$ and the Sobolev norms of $\partial^k_t(\bar{\rho}, U,V,U_x), k\leq k_0$.
Setting
\begin{align}
\label{4.4}
u^{a1}(t,x,\eta)=U(t,x)+u^a(t,x,\eta), \quad v^{a1}(t,x,\eta)=V(t,x)\eta+v^a(t,x,\eta),
\end{align}
then, $(u^{a1}, v^{a1})$ is an approximate solution to the problem \eqref{1.1} satisfying compatibility conditions up to order $k_0$ and initial data.

\subsubsection{Improving decay in $\eta$}

Note that the approximate solution  $(u^{a1}, v^{a1})$ satisfies
\begin{equation}\label{far}
\lim\limits_{\eta\to+\infty} u^{a1}(t,x,\eta)=U(t,x),\end{equation}
and the divergence constraint
\begin{equation}\label{div}
\partial_x(\bar{\rho}u^{a1})+\partial_{\eta}(\bar{\rho}v^{a1})=-\bar{\rho}_t,
\end{equation}
for all $t\ge 0$. However, the error 
$$f^{a1}=(\partial_t+u^{a1}\partial_x
+v^{a1}\partial_\eta-\frac{1}{\bar{\rho}(t,x)}\partial_{\eta}^2)u^{a1}+P_x$$
does not have enough
 decay compared with  $\partial_\eta u^{a1}$  as $\eta\to +\infty$. 
Since this property is essential
 for the convergence of the Nash-Moser-H\"ormander iteration scheme of the nonlinear problem given  in next section, we need to modify  the approximate solution  $(u^{a1}, v^{a1})$ as follows.

From \eqref{1.1},  $\partial_\eta u$ satisfies
$$\begin{cases}
(\partial_t+u\partial_x
+v\partial_\eta-\frac{1}{\bar{\rho}(t,x)}\partial_{\eta}^2)u_\eta+(u_x+v_\eta)u_\eta=0,\\
 u_{\eta}|_{\eta=0}=\bar{\rho}P_x.
\end{cases}
$$
This motivates us to consider the  following initial-boundary value problem for a linear degenerate parabolic equation:
\begin{align}
\label{APP02}
\left\{
\begin{array}{ll}
\phi_t+(u^{a1}\phi)_x+(v^{a1}\phi)_{\eta}-\frac{1}{\bar{\rho}}\partial^2_{\eta}\phi=0,\\
\partial_{\eta}\phi|_{\eta=0}=(\bar{\rho}P_x)(t,x),\\
\phi|_{t=0}=(\partial_\eta u_{0})(x,\eta).
\end{array}
\right.
\end{align}

Suppose that the solution $\phi$ of \eqref{APP02} is obtained.
 Define an approximate solution $(u^{a2}, v^{a2})$ as 
\begin{align}
\label{APP03}
u^{a2}=U(t,x)-\int_{\eta}^{\infty}\phi(t,x,\tilde{\eta})d\tilde{\eta},\quad v^{a2}=V\eta+\frac{1}{\bar{\rho}}\int_0^{\eta}(\bar{\rho}\int_{\tilde{\eta}}^{\infty}\phi(t,x,s)ds)_xd\tilde{\eta}.
\end{align}
It is straightforward to verify that the compatibility conditions of \eqref{1.1}, the far-field condition \eqref{far} 
and the divergence constraint \eqref{div} still hold for  $(u^{a2}, v^{a2})$. Moreover, it satisfes the equation,
\begin{align}
\label{APP04}
u^{a2}_t+u^{a2}u^{a2}_x+v^{a2}u^{a2}_{\eta}+P_x-\frac{1}{\bar{\rho}}\partial^2_{\eta}u^{a2}=f^0,
\end{align}
where
$$
f^0=-\int_{\eta}^{\infty}\left(\left[U-u^{a1}-\int_{\tilde{\eta}}^{\infty}\phi(t,x,s)ds\right]\phi\right)_xd\tilde{\eta}-(v^{a2}-v^{a1})\phi,
$$
which will be shown to decay
 faster than $\partial_\eta  u^{a2}$.\\

%To study properties of the solution of (\ref{APP02}), we introduce the following Sobolev norm:
%\begin{align}
%\label{APP0}
%\|\phi\|_{H^{k_0}_{\gamma}}=\sum_{|\alpha|\leq k_0}\|(1+\eta)^{\gamma+\alpha_2}D^{\alpha}\phi\|%_{L^2(\mathbb{T}\times\mathbb{R}^+)}
%\end{align}
%where $D^{\alpha}=\partial_x^{\alpha_1}\partial_{\eta}^{\alpha_2}, \alpha=(\alpha_1,\alpha_2)$. 
From the boundedness of $u^a$ given in (\ref{APP06}) and
some elementary weighted energy estimates on
the  solution to (\ref{APP02}), we have
\begin{pro}\label{prop4.1}
Under the Main Assumptions (H) on the initial data, there exists a unique solution $\phi(t,x,\eta)$ to (\ref{APP02}). Moreover, there is $T>0$ such that $\phi$ satisfies 
\begin{align}
\label{APP07}
\begin{cases}
\max_{0\le t\le T}\|\phi(t)\|_{H^{2k_0}_{\gamma}}\leq C_1, \quad  \phi(t,x,\eta)\geq \frac{C_2}{(1+\eta)^{\gamma+2}},\quad \forall (t,x,\eta)\in [0,T]\times \mathbb{T}\times\mathbb{R}^+
,\\[3mm]
\|(1+\eta)^{\gamma+2+\alpha_2}D^{\alpha}\phi\|_{L^{\infty}([0,T]\times \mathbb{T}\times\mathbb{R}^+)}\leq C_3,\quad |\alpha|\leq k_0,
\end{cases}
\end{align}
where
$$\|\phi(t)\|_{H^{2k_0}_{\gamma}}=\sum_{|\alpha|\leq 2k_0}\|(1+\eta)^{\gamma+\alpha_2}D^{\alpha}\phi(t)\|_{L^2(\mathbb{T}\times\mathbb{R}^+)},$$
with $D^{\alpha}=\partial_{t,x}^{\alpha_1}\partial_{\eta}^{\alpha_2}, \alpha=(\alpha_1,\alpha_2)$. 

\end{pro}
\begin{pf}
The proof is divided in two steps.

(1) Applying the operator $D^{\alpha}=\partial_{t,x}^{\alpha_1}\partial_{\eta}^{\alpha_2}$ to the equation $(\ref{APP02})$, multiplying the resulting equation by $(1+\eta)^{2\gamma+2\alpha_2}D^{\alpha}\phi$ and integrating it over $\mathbb{T}\times\mathbb{R}^+$, we obtain
\begin{align}\label{4.11}
\frac12\frac{d}{dt}\|(1+\eta)^{\gamma+\alpha_2}D^{\alpha}\phi(t)\|^2_{L^2}+\|\frac{1}{\sqrt{\bar{\rho}}}(1+\eta)^{\gamma+\alpha_2}D^{\alpha}\phi_{\eta}(t)\|^2_{L^2}=\sum_{i=1}^7 I_i,
\end{align}
where
$$
I_1=\int_{\mathbb{T}\times\mathbb{R}^+}(u^{a1}_x+v^{a1}_{\eta})(1+\eta)^{2\gamma+2\alpha_2}(D^{\alpha}\phi)^2 dxd\eta,$$
$$
I_2=\int_{\mathbb{T}\times\mathbb{R}^+}(u^{a1}D^{\alpha}\phi_x+v^{a1}D^{\alpha}\phi_{\eta})(1+\eta)^{2\gamma+2\alpha_2}D^{\alpha}\phi dxd\eta,$$
$$
I_3=2(\gamma+\alpha_2)\int_{\mathbb{T}\times\mathbb{R}^+}\frac{1}{\bar{\rho}}(1+\eta)^{2\gamma+2\alpha_2-1}(D^{\alpha}\phi_{\eta})(D^{\alpha}\phi) dxd\eta,$$
$$
I_4=\sum_{0<\beta\le \alpha}C_{\alpha}^{\beta}\int_{\mathbb{T}\times\mathbb{R}^+} (1+\eta)^{2\gamma+2\alpha_2} D^{\beta}(u^{a1}_x+v^{a1}_{\eta})(D^{\alpha-\beta}\phi)
(D^{\alpha}\phi) dxd\eta,$$
$$
I_5=\sum_{0<\beta\le \alpha}C_{\alpha}^{\beta}\int_{\mathbb{T}\times\mathbb{R}^+}(1+\eta)^{2\gamma+2\alpha_2}(D^{\beta}u^{a1}D^{\alpha-\beta}\phi_x
+D^{\beta}v^{a1}D^{\alpha-\beta}\phi_{\eta})
D^{\alpha}\phi dxd\eta, 
$$
$$
I_6=\sum_{0<\beta\le \alpha, \beta_2=0}C_{\alpha}^{\beta}\int_{\mathbb{T}\times\mathbb{R}^+}(1+\eta)^{2\gamma+2\alpha_2}(D^{\beta}\frac{1}{\bar{\rho}})(\partial_{\eta}^2
D^{\alpha-\beta}\phi) (D^{\alpha}\phi) dxd\eta,
$$
and
$$I_7=\int_{\mathbb{T}}\frac{1}{\bar{\rho}}(\partial_{\eta}D^{\alpha}\phi D^{\alpha}\phi)|_{\eta=0}dx.$$

It is straightforward  to show
\begin{align*}
|I_1|
\leq&\|u^{a1}_x,v^{a1}_{\eta}\|_{L^{\infty}}\|(1+\eta)^{\gamma+\alpha_2}D^{\alpha}\phi\|^2_{L^2},
\end{align*}
\begin{align*}
|I_2|=&|\int_{\mathbb{T}\times\mathbb{R}^+}(u^{a1}D^{\alpha}\phi_x+v^{a1}D^{\alpha}\phi_{\eta})(1+\eta)^{2\gamma+2\alpha_2}D^{\alpha}\phi dxd\eta|\\
\lesssim& (\|u^{a1}_x\|_{L^{\infty}}+\|v^{a1}_{\eta}\|_{L^{\infty}}+\|v^{a1}_{\eta}/(1+\eta)\|_{L^{\infty}})\|(1+\eta)^{\gamma+\alpha_2}D^{\alpha}\phi\|^2_{L^2},
\end{align*}
by integration by parts, and using
$v^{a1}|_{\eta=0}=0$, and
\begin{align*}
|I_3|=&2(\gamma+\alpha_2)|\int_{\mathbb{T}\times\mathbb{R}^+}\frac{1}{\bar{\rho}}D^{\alpha}\phi_{\eta}(1+\eta)^{2\gamma+2\alpha_2-1}D^{\alpha}\phi dxd\eta|\\
\leq&\frac18\|\frac{1}{\sqrt{\bar{\rho}}}(1+\eta)^{\gamma+\alpha_2}D^{\alpha}\phi_{\eta}\|^2_{L^2}+C\|\frac{1}{\bar{\rho}}\|_{L^{\infty}}\|(1+\eta)^{\gamma+\alpha_2}
D^{\alpha}\phi\|^2_{L^2}.
\end{align*}

On the other hand, we have
\begin{align*}
|I_4|
=&|\sum_{0<\beta\le \alpha}C_{\alpha}^{\beta}\int_{\mathbb{T}\times\mathbb{R}^+}
(1+\eta)^{\beta_2}D^{\beta}(\frac{-\bar{\rho}_t-\bar{\rho}_xu^{a1}}{\bar{\rho}})
(1+\eta)^{\gamma+\alpha_2-\beta_2}D^{\alpha-\beta}\phi
(1+\eta)^{\gamma+\alpha_2}
D^{\alpha}\phi dxd\eta|\\
\leq&C\|\sum_{\beta}D^{\beta}(\frac{\bar{\rho}_t+\bar{\rho}_xU}{\bar{\rho}})\|_{L^{\infty}}\|\phi\|^2_{H^{2k_0}_{\gamma}}+|\tilde{I}_4|,
\end{align*}
by noting that
$\beta_2=0$ for the operator $D^\beta$ acting on $\frac{\bar{\rho}_t+\bar{\rho}_xU}{\bar{\rho}}$. Hence,
$$
|\tilde{I}_4|\leq C\|(1+\eta)^{\beta_2}\left(D^{\beta}\frac{\bar{\rho}_x(u^{a1}-U)}{\bar{\rho}}\right)\|_{L^{\infty}}\|\phi\|^2_{H^{2k_0}_{\gamma}},\ \  \hbox{for}\ \ |\beta|\leq k_0,
$$
and
\begin{align*}
|\tilde{I}_4|\leq& C\|(1+\eta)^{\beta_2} (D^{\beta}\frac{\bar{\rho}_x(u^{a1}-U)}{\bar{\rho}})\|_{L^2}\|(1+\eta)^{\gamma+\alpha_2-\beta_2}D^{\alpha-\beta}\phi\|_{L^{\infty}}
\|(1+\eta)^{\gamma+\alpha_2} D^{\alpha}\phi\|_{L^2}\\
\leq& C\|(1+\eta)^{\beta_2}\left(D^{\beta}\frac{\bar{\rho}_x(u^{a1}-U)}{\bar{\rho}}\right)\|_{L^2}\|\phi\|^2_{H^{2k_0}_{\gamma}},\ \  \hbox{for}\ \ |\beta|> k_0,
\end{align*}
by using the weighted Sobolev embedding. Similarly, one has
\begin{align*}
|I_5|=&|\sum_{0<\beta\le \alpha}\int_{\mathbb{T}\times\mathbb{R}^+}C_{\alpha}^{\beta}(D^{\beta}u^{a1}D^{\alpha-\beta}\phi_x
+D^{\beta}v^{a1}D^{\alpha-\beta}\phi_{\eta})(1+\eta)^{2\gamma+2\alpha_2}
D^{\alpha}\phi dxd\eta|\\
\leq&|I_{5}^1|+|I^2_5|,
\end{align*}
where
$$
|I_{5}^1|\leq 
\begin{cases}
C(\|D^{\beta}(u^{a1}-U)(1+\eta)^{\beta_2}\|_{L^{\infty}}+\|D^{\beta}U\|_{L^{\infty}})\|\phi\|^2_{H^{2k_0}_{\gamma}},\ \ \hbox{for}\ \ 1\leq|\beta|\leq k_0, \\[2mm]
 C(\|D^{\beta}(u^{a1}-U)(1+\eta)^{\beta_2}\|_{L^2}+\|D^{\beta}U\|_{L^{\infty}})\|\phi\|^2_{H^{2k_0}_{\gamma}},\ \ \hbox{for}\ \ |\beta|> k_0,
\end{cases}
$$
and
$$
|I_5^2|\leq 
\begin{cases}
C\|D^{\beta}v^{a1}(1+\eta)^{\beta_2-1}\|_{L^{\infty}}\|\phi\|^2_{H^{2k_0}_{\gamma}},\ \ \hbox{for}\ \ 1\leq|\beta|\leq k_0, \\[2mm]
C\|D^{\beta}v^{a1}(1+\eta)^{\beta_2-1}\|_{L^2}\|\phi\|^2_{H^{2k_0}_{\gamma}},\ \ \hbox{for}\ \ \beta_2\neq0,\ |\beta|> k_0;
\end{cases}
$$
\begin{align*}
|I_5^2|\leq& C\|D^{\beta}v^{a1}(1+\eta)^{-1}\|_{L^{\infty}_{\eta}(L^2_x)}\|D^{\alpha-\beta}\phi_{\eta}(1+\eta)^{\gamma+\alpha_2+1}\|_{L^2_{\eta}(L^{\infty}_x)}
\|D^{\alpha}\phi(1+\eta)^{\gamma+\alpha_2}\|_{L^2},\\
\leq&C\|D^{\beta}v^{a1}(1+\eta)^{-1}\|_{L^{\infty}_{\eta}(L^2_x)}\|\phi\|^2_{H^{2k_0}_{\gamma}}\ \ \hbox{for}\ \ \beta_2=0,\ |\beta|> k_0.
\end{align*}

For the term $I_6$, by integration by parts, we have
\begin{align*}
|I_6|
=&\sum_{0<\beta\le \alpha, \beta_2=0}\left(|
\int_{\mathbb{T}\times\mathbb{R}^+}C_{\alpha}^{\beta}(2\gamma+2\alpha_2)(1+\eta)^{2\gamma+2\alpha_2-1}
(D^{\beta}\frac{1}{\bar{\rho}})(\partial_{\eta}D^{\alpha-\beta}\phi)
(D^{\alpha}\phi) dxd\eta \right.\\
&
+\int_{\mathbb{T}\times\mathbb{R}^+}C_{\alpha}^{\beta}(1+\eta)^{2\gamma+2\alpha_2}(D^{\beta}\frac{1}{\bar{\rho}})(\partial_{\eta}
D^{\alpha-\beta}\phi)(D^{\alpha}\phi_{\eta})dxd\eta|\\
& \left.
-|\int_{\mathbb{T}}C_{\alpha}^{\beta}(D^{\beta}\frac{1}{\bar{\rho}})(\partial_{\eta}D^{\alpha-\beta}\phi) (D^{\alpha}\phi)|_{\eta=0}dx|\right)\\
\leq&
\sum_{0<\beta\le \alpha, \beta_2=0}\left(
|\int_{\mathbb{T}}C_{\alpha}^{\beta}(D^{\beta}\frac{1}{\bar{\rho}})(\partial_{\eta}D^{\alpha-\beta}\phi) (D^{\alpha}\phi)|_{\eta=0}dx|+C\|(1+\eta)^{\gamma+\alpha_2}D^{\alpha-\beta}\phi_{\eta}\|_{L^2}^2\right)\\
&
+\frac{1}{8}\|\frac{1}{\sqrt{\bar{\rho}}}(1+\eta)^{\gamma+\alpha_2}D^{\alpha}\phi_{\eta}\|_{L^2}^2.
\end{align*}

It remains to handle the boundary integration terms on the right hand side of
the above estimate and $I_7$.
For illustration,  we only estimate $I_7$.\\
Firstly, noticing that
\begin{align*}
\partial_{\eta}\phi|_{\eta=0}=\bar{\rho}P_x,
\end{align*}
and applying the operator $\partial_{\eta}$ on the equation (\ref{APP02}), we obtain
\begin{align*}
\frac{1}{\bar{\rho}}\partial^3_{\eta}\phi=(\phi_{\eta})_t+(u^{a1}_x\phi)_{\eta}+(u^{a1}\phi_x)_{\eta}+(v^{a1}\phi)_{\eta\eta}.
\end{align*}
Taking  this equation on the boundary $\{\eta=0\}$  and using the boundary condition, we get
\begin{align*}
\frac{1}{\bar{\rho}}\partial^3_{\eta}\phi|_{\eta=0}=(\bar{\rho}P_x)_t+u^{a1}_x\bar{\rho}P_x+u^{a1}_{x\eta}\phi|_{\eta=0}+u^{a1}(\bar{\rho}P_x)_x+
u^{a1}_{\eta}\phi_{x}|_{\eta=0}+2v^{a1}_{\eta}\bar{\rho}P_x+v^{a1}_{\eta\eta}\phi|_{\eta=0}.
\end{align*}
By induction, for positive integer  $k$, we have
\begin{align*}
\frac{1}{\bar{\rho}}\partial^{2k+1}_{\eta}\phi|_{\eta=0}=&(\partial_{\eta}^{2k-1}\phi)_t+\partial_{\eta}^{2k-1}(u^{a1}_x\phi)
+\partial_{\eta}^{2k-1}(u^{a1}\phi_x)
+\partial_{\eta}^{2k}(v^{a1}\phi)\\
=&[F(D^{\alpha}_{|\alpha|\leq 2k-2}u^{a1},D^{\beta}_{|\beta|\leq 2k-2}v^{a1},D^{\gamma}_{|\gamma|\leq 2k-3}\phi, D^{\pi}(\bar{\rho}, P_x))]_t
+\sum_{i=1}^{2k-1}C_{2k-1}^i\partial_{\eta}^iu^{a1}_x\partial_{\eta}^{2k-1-i}\phi\\
&+\sum_{j=1}^{2k-1}C_{2k-1}^j\partial_{\eta}^iu^{a1}\partial_{\eta}^{2k-1-j}\phi_x
+\sum_{s=1}^{2k}C_{2k}^s\partial_{\eta}^sv^{a1}\partial_{\eta}^{2k-s}\phi\\
=&G(D^{\alpha}_{|\alpha|\leq 2k}u^{a1},D^{\beta}_{|\beta|\leq 2k}v^{a1},D^{\gamma}_{|\gamma|\leq 2k-1}\phi, D^{\pi}(\bar{\rho}, P_x))),
\end{align*}
where $F, G$ are polynomial functions.  Hence,  the normal derivative of $\phi$ can be reduced by two order using
 the boundary condition and the equation (\ref{APP02}).  Therefore,
 we can use the trace estimate to control the boundary integral.

\iffalse

And by the standard Sobolev embedding inequality $\|\partial_{\eta}^{2k-1}\phi\|_{L^{\infty}(\mathbb{T}\times\mathbb{R}^+)}\leq C\|\partial_{\eta}^{2k-1}\phi\|_{H^{2}(\mathbb{T}\times\mathbb{R}^+)}$. Similarly, the $L^{\infty}$ norms of other terms $D^{\alpha}u^{a1}, D^{\beta}v^{a1}, D^{\pi}(\bar{\rho}, P_x)$ can be bounded by their corresponding Sobolev norms. By the definitions of $u^{a1}, v^{a1}$ and the Main Assumptions (H),

\fi

Thus, by summarizing the above estimates, and taking summation over $|\alpha|\leq 2k_0$ for \eqref{4.11}, it follows
\begin{align*}
\frac{d}{dt}\|\phi(t)\|^2_{H^{2k_0}_{\gamma}}+\sum_{|\alpha|\leq 2k_0}\|\frac{1}{\sqrt{\bar{\rho}}}(1+\eta)^{\gamma+\alpha_2}D^{\alpha}\phi_{\eta}(t)\|^2_{L^2}\leq C\|\phi(t)\|^2_{H^{2k_0}_{\gamma}}.
\end{align*}
which implies the first boundedness estimate given in (\ref{APP07}) by using Gronwall inequality.\\

(2)
Next, we apply the maximal principle to prove the second estimate
 given in (\ref{APP07}).

From \eqref{APP02},  $y(t,x,\eta)\triangleq(1+\eta)^{\gamma+2}\phi$ satisfies the following degenerate parabolic equation,
\begin{align*}
y_t+(u^{a1}_x+v^{a1}_{\eta}-\frac{v^{a1}(2+\gamma)}{1+\eta}-\frac{(\gamma+2)(\gamma+3)}{\bar{\rho}(1+\eta)^2})y+u^{a1}y_x
+(v^{a1}+\frac{2(\gamma+2)}{\bar{\rho}(1+\eta)})y_{\eta}-\frac{1}{\bar{\rho}}\partial^2_{\eta}y=0.
\end{align*}
By the maximal principle (see also Lemma E.2 in \cite{MW1}), we have
\begin{align*}
\min_{\mathbb{T}\times\mathbb{R}^+}y(t)\geq(1-\lambda te^{\lambda t})k(t),
\end{align*}
with
$$
k(t)=\min\{\min_{\mathbb{T}\times\mathbb{R}^+}y|_{t=0}, \min_{[0,t]\times\mathbb{T}}y|_{\eta=0}\},
$$
for a fixed $\lambda\geq \|(u^{a1}_x+v^{a1}_{\eta}-\frac{v^{a1}(2+\gamma)}{1+\eta}-\frac{(\gamma+2)(\gamma+3)}{\bar{\rho}(1+\eta)^2})\|_{L^{\infty}}$.

It follows from the Main Assumptions (H2) on the initial data that
$$
\min_{\mathbb{T}\times\mathbb{R}^+}y(0)\geq \sigma_0>0.
$$
It suffices to derive the lower bound on $\min_{[0,t]\times\mathbb{T}}y|_{\eta=0}$. Notice that $y|_{\eta=0}=\phi|_{\eta=0}$, the first boundedness estimate of  (\ref{APP07}) and the Sobolev inequality give 
$$
\|\phi_t|_{\eta=0}\|_{L^{\infty}}\leq C.
$$
Consequently,
\begin{align*}
\phi(t)|_{\eta=0}\geq \phi(0,x,\eta)|_{\eta=0}-Ct\geq \sigma_0-Ct.
\end{align*}
Thus,  we have the lower bound given in the second 
estimate in (\ref{APP07}) provided that $t$ is suitably small. The 
third estimate in  $(\ref{APP07})$ can also be proved by the  maximal principle similarly (also refer to Lemma E.1 in \cite{MW1}).   Then the
 proof of this proposition is completed.
\end{pf}\\

\iffalse

The second level approximate solution $(u^{a2}, v^{a2})$ is defined as follows.
\begin{align}
\label{APP03}
u^{a2}=U(t,x)-\int_{\eta}^{\infty}\phi(t,x,\tilde{\eta})d\tilde{\eta},\quad v^{a2}=V\eta+\frac{1}{\bar{\rho}}\int_0^{\eta}(\bar{\rho}\int_{\tilde{\eta}}^{\infty}\phi(t,x,s)ds)_xd\tilde{\eta}.
\end{align}
It is direct to check that the definitions in (\ref{APP03}) are well-defined due to the estimates (\ref{APP07}). Moreover, the further approximate solution $(u^{a2}, v^{a2})$ satisfies the following system.
\begin{align}
\label{APP04}
\left\{
\begin{array}{ll}
u^{a2}_t+u^{a2}u^{a2}_x+v^{a2}u^{a2}_{\eta}+P_x-\frac{1}{\bar{\rho}}\partial^2_{\eta}u^{a2}=f^0;\\
\partial_{\eta}(\bar{\rho}v^{a2})+\partial_x(\bar{\rho}u^{a2})=-\bar{\rho}_t,\\
u^{a2}(t=0,x,\eta)=u_0(x,\eta)
\end{array}
\right.
\end{align}
here
$$
f^0=-\int_{\eta}^{\infty}\left(\left[U-u^{a1}-\int_{\tilde{\eta}}^{\infty}\phi(t,x,s)ds\right]\phi\right)_xd\tilde{\eta}-(v^{a2}-v^{a1})\phi
$$
It follows from the definition (\ref{APP03}) that
$u^{a2}_{\eta}=\phi$. Here we should remark that the error term $f^0$ is much more accurate than that of the first approximate solution $(u^{a1}, v^{a1})$ in the sense that it has higher order decay rates with respect to the fast variable $\eta$, because the error terms $f^0$ come from the nonlinear terms and every linear term has decay properties of (\ref{APP07}).\\

\fi

\subsubsection{Boundary condition}

It is noted that the  approximate solution $u^{a2}$ does not satisfy the original boundary condition, that is, $u^{a2}|_{\eta=0}\neq 0$. For this,
set $$\zeta(t,x)\triangleq u^{a2}|_{\eta=0}=U(t,x)-\int_0^{\infty}\phi(t,x,\eta)d\eta.$$
 $\zeta(t,x)$ is uniformly continuous and bounded due to (\ref{APP07}). By the compatibility condition of the initial data, we have $\zeta(0,x)=0$. Consequently, $|\zeta(t,x)|\leq \varepsilon_0,\ t\in[0,t_0]$, with $\varepsilon_0\rightarrow 0$ as $t_0$ tends to zero.\\
In addition, there exists a smooth monotone decreasing function $\psi(\eta)\subseteq[0,1], \eta\geq 0$ such that $\hbox{supp}\psi\subseteq [0,1]$,  $\psi(0)=1$ and $|\psi'(\eta)|<C$.
Note that there exists a positive constant $a_0$ such that $\phi(t,x,\eta)>a_0, \eta\in[0,1]$. 

Now,  define 
\begin{align}
\label{APP09}
u^{a3}=u^{a2}-\zeta(t,x)\psi(\eta), \quad v^{a3}=-\frac{1}{\bar{\rho}}\int_0^{\eta}(\bar{\rho}(u^{a2}(t,x,\tilde{\eta})-U(x,t)-\zeta(t,x)\psi(\tilde{\eta}))_xd\tilde{\eta}.
\end{align}
It is direct to check that $u^{a3}(t,x,\eta)|_{\eta=0}=0$, and
\begin{align}
\label{APP20}
u^{a3}_{\eta}(t,x,\eta)=\phi(t,x,\eta)-\zeta(t,x)\psi'(\eta)>\frac{\phi}{2}>0,
\end{align}
provided that $t\in[0,t_0]$ with $ t_0$ being
 suitably small. And the profile $(u^{a3}, v^{a3})$ satisfies
\begin{align}
\label{APP10}
\left\{
\begin{array}{ll}
u^{a3}_t+u^{a3}u^{a3}_x+v^{a3}u^{a3}_{\eta}+P_x-\frac{1}{\bar{\rho}}\partial^2_{\eta}u^{a3}=f^a,\\
\partial_{\eta}(\bar{\rho}v^{a3})+\partial_x(\bar{\rho}u^{a3})=-\bar{\rho}_t,\\
u^{a3}(0,x,\eta)=(\partial_\eta u_0)(x,\eta),
\end{array}
\right.
\end{align}
with
$f^a=f^0-\bar{f}^0$, where
\begin{align*}
\bar{f}^0=\zeta_t\psi+\zeta\psi u^{a2}_x+u^{a2}\zeta_x\psi-\zeta\zeta_x\psi^2
+v^{a2}\zeta\psi'-\frac{u^{a2}_{\eta}(\bar{\rho}\zeta)_x}{\bar{\rho}}\int_0^{\eta}\psi(\tilde{\eta})d\tilde{\eta}
+\frac{\zeta\psi'(\bar{\rho}\zeta)_x}{\bar{\rho}}\int_0^{\eta}\psi(\tilde{\eta})d\tilde{\eta}-\frac{\psi^{''}\zeta}{\bar{\rho}}.
\end{align*}
\begin{re}
The  approach of constructing the zero-th approximate solution to (\ref{1.1}) introduced above can be applied
to  the incompressible Prandtl equations. 
\end{re}

\subsection{The Nash-Moser-H\"omander iteration scheme}

We now construct the approximate solution sequence of \eqref{1.1} by using the Nash-Moser-H\"omander Iteration Scheme. The procedure mainly follows the  one  given in \cite{AWXY}.
Thus,  we will only present the main steps.
  
Denote the linearized operator $\mathcal{P}'$ around
 $(\hat{\omega}, \hat{q})$ of (\ref{1.1}) by 
$$
\mathcal{P}'_{(\hat{\omega}, \hat{q})}(\omega,
q)=\partial_t\omega+\hat{\omega}\omega_x+\hat{q}\omega_{\eta}+\omega\hat{\omega}_x+q\hat{\omega}_{\eta}-\frac{1}{\bar{\rho}}\partial_{\eta}^2\omega.
$$

Suppose that the approximate solutions $(u^k, v^k) $ of \eqref{1.1} have been
constructed for all $k\le n$, with $u^0=u^{a3}$ and $v^0=v^{a3}$ being defined in Subsection 4.1.3, we construct the $(n+1)-$th approximate solution
$(u^{n+1}, v^{n+1})$ as follows:
\begin{align}
\label{4.1}
u^{n+1}=u^{n}+\delta u^n=u^{a3}+\tilde{u}^n+\delta u^n, \quad
v^{n+1}=v^n+\delta v^n=v^{a3}+\tilde{v}^n+\delta v^n,
\end{align}
where the increment $(\delta u^n, \delta v^n)$ is the solution
to the following initial-boundary value problem
\begin{align}
\label{4.2} \left\{
\begin{array}{ll}
\mathcal{P}'_{(u_{\theta_n}^n, v_{\theta_n}^n)}(\delta u^n, \delta v^n)=f^n,\\
\partial_{\eta}(\bar{\rho}(t,x)\delta v^n)+\partial_x(\bar{\rho}(t,x)\delta u^n)=0,\\
\delta u^n|_{\eta=0}=\delta v^n|_{\eta=0}=0,\quad
\lim_{\eta\rightarrow+\infty}\delta u^n=0,\\
\delta u^n|_{t=0}=0.
\end{array}
\right.
\end{align}
Here,  $u_{\theta_n}^n=u^{a3}+S_{\theta_n}\tilde{u}^n$ and
$v_{\theta_n}^n=v^{a3}+S_{\theta_n}\tilde{v}^n$ with
$\theta_n=\sqrt{\theta_0^2+n}$ for any $n\geq 1$ and a large fixed
constant $\theta_0$. The smoothing operator $S_{\theta}$ is defined
by
$$
(S_{\theta}f)(t,x,\eta)=\iiint
j_{\theta}(\tau)j_{\theta}(\xi)j_{\theta}(\mu)\tilde{f}(t-\tau+\theta^{-1},
x-\xi,\eta-\mu+\theta^{-1})d\tau d\xi d\mu,
$$
for a function  $f$ defined on
$\Omega=[0,+\infty[\times\mathbb{T}_x\times\mathbb{R}_{\eta}^+$ with
$\tilde{f}$ being the zero extension of $f$ to $\mathbb{R}^3$, and the
mollifier $j_{\theta}(\tau)=\theta j(\theta\tau)$ with $ j\in
C_0^{\infty}(\mathbb{R})$  being 
a non-negative function satisfying $\hbox{Supp}j\subseteq[-1,1]$ and
$\|j\|_{L^1}=1$.

In order to show that the approximate solution $(u^n, v^n)$
converges to the solution of the nonlinear problem (\ref{1.1}), we need
to define
the source term $f^n$ properly for the problem (\ref{4.2}).

To do this, denoting the nonlinear operator on the left hand side
of (\ref{1.1}) by $\mathcal{P}(\omega, q)$, obviously, the following
identity holds:
\begin{align}
\label{4.3} \mathcal{P}(u^{n+1},
v^{n+1})-\mathcal{P}(u^{n},
v^{n})=\mathcal{P}'_{(u_{\theta_n}^n,
v_{\theta_n}^n)}(\delta u^n,\delta v^n)+e_n,
\end{align}
where
$$
e_n=e_n^1+e_n^2,
$$
with $e^1_n$ being the error term from the Newton iteration,
\begin{align} \label{e-1}
e_n^1=&\mathcal{P}(u^n+\delta u^n, v^n+\delta
v^n)-\mathcal{P}(u^n, v^n)-\mathcal{P}'_{(u^n,
v^n)}(\delta u^n,\delta v^n)\nonumber\\
=&\delta u^n\partial_x(\delta u^n)+\delta
v^n\partial_{\eta}(\delta u^n),
\end{align}
and $e_n^2$ being the error from mollifying the coefficients,
\begin{align}\label{e-2}
e^2_n=&\mathcal{P}'_{(u^n, v^n)}(\delta u^n, \delta
v^n)-\mathcal{P}'_{(u^n_{\theta_n},
v_{\theta_n}^n)}(\delta u^n, \delta v^n)\nonumber\\
=&((1-S_{\theta_n})(u^n-u^{a3}))\partial_x(\delta u^n)+\delta u^n\partial_x((1-S_{\theta_n})(u^n-u^{a3}))\nonumber\\
&+((1-S_{\theta_n})(v^n-v^{a3}))\partial_{\eta}\delta u^n+\delta
v^n\partial_{\eta}((1-S_{\theta_n})(u^n-u^{a3})).
\end{align}
Taking summation of  (\ref{4.3}) over all $n\in\mathbb{N}$ leads to
\begin{align}
\label{4.4} \mathcal{P}(u^{n+1},
v^{n+1})=\sum_{j=0}^n(\mathcal{P}'_{(u_{\theta_j}^j,
v^j_{\theta_j})}(\delta u^j, \delta v^j)+e_j)+f^{a},
\end{align}
with $f^a=\mathcal{P}(u^{a3},v^{a3})$.

It is obvious that if  the approximate solution $(u^n,
v^n)$ converges to the solution to (\ref{1.1}), then the right hand side of  (\ref{4.4}) must converge to zero as $n$ tends to $+\infty$. In
this way, it is convenient to require that $(\delta u^n, \delta v^n)\
(n\geq 0)$ satisfies the equation,
$$
\mathcal{P}'_{(u_{\theta_n}^n, v_{\theta_n}^n)}(\delta u^n,
\delta v^n)=f^n,
$$
where $f^n$ is defined by
\begin{align}
\label{EQU}
\sum_{j=0}^nf^j=-S_{\theta_n}(\sum_{j=0}^{n-1}e_j)-S_{\theta_n}f^a,
\end{align}
inductively, that is,
\begin{align}
\label{4.5} \left\{
\begin{array}{ll}
f^0=-S_{\theta_0}f^a,\quad
f^1=(S_{\theta_0}-S_{\theta_1})f^a+S_{\theta_0}f^a,\\
f^n=(S_{\theta_{n-1}}-S_{\theta_n})(\sum_{j=0}^{n-2}e_j)-S_{\theta_{n}}e_{n-1}+(S_{\theta_{n-1}}-S_{\theta_n})f^a,
\quad \forall n\geq 2,
\end{array}
\right.
\end{align}
with $f^a$  given in \eqref{APP10}.

We now give some properties of the smoothing
operator in the following lemma, which also can be found in  Section 4.1 of \cite{AWXY}.
\begin{lem}
The smoothing operator $\{S_{\theta}\}_{\theta>0}:
\mathcal{A}_l^0(\Omega)\rightarrow \cap_{s\geq
0}\mathcal{A}_l^s(\Omega)$, satisfies the following estimates:
\begin{align}
\label{4.6} \left\{
\begin{array}{ll}
\|S_{\theta}v\|_{\mathcal{A}_l^s}\leq
C_j\theta^{(s-\alpha)_+}\|v\|_{\mathcal{A}_l^{\alpha}}, \quad \hbox{for
all}\ s, \alpha\geq 0,\\
\|(1-S_{\theta})v\|_{\mathcal{A}_l^s}\leq
C_j\theta^{(s-\alpha)}\|v\|_{\mathcal{A}_l^{\alpha}}, \quad \hbox{for
all}\ 0\leq s\leq \alpha,
\end{array}
\right.
\end{align}
and
\begin{align}
\label{4.7}
\|(S_{\theta_n}-S_{\theta_{n-1}})v\|_{\mathcal{A}^s_l}\leq
C_j\triangle\theta_{n}^{s-\alpha}\|v\|_{\mathcal{A}_l^{\alpha}},
\hbox{for all}\ s, \alpha\geq 0,
\end{align}
where $\triangle\theta_n=\theta_{n+1}-\theta_n$, and the constant
$C_j$ depends only on the mollifier function $j\in
C_0^{\infty}(\mathbb{R})$.
\end{lem}

\iffalse

And we also have the following commutator estimates
\begin{lem}
For any proper function $v$, it follows that
\begin{align}
\label{4.8} \|[\frac{1}{\partial_{\eta}u^{a3}},
S_{\theta}](\partial_{\eta}v)\|_{\mathcal{A}_l^k}\leq
C_k\|\frac{v}{\partial_{\eta}u^{a3}}\|_{\mathcal{A}_l^k}
\end{align}
and
\begin{align}
\label{4.9} \|\partial_{\eta}[\frac{1}{\partial_{\eta}u^{a3}},
\partial_{\eta}S_{\theta}]v\|_{\mathcal{A}_l^k}\leq
C_k\theta\|\frac{v}{\partial_{\eta}u^{a3}}\|_{\mathcal{A}_l^k}
\end{align}
here the constant $C_k$ depends on $W$. Moreover, similar estimates
on the norms $\|\cdot\|_{\mathcal{B}^{k_1, k_2}_l},
\|\cdot\|_{\mathcal{C}^{k}_l}$ and $\|\cdot\|_{\mathcal{D}^{k}_l}$.
\end{lem}
\begin{pf}
The proof of this lemma is similar as Lemma 4.1 in \cite{AWXY}. The only difference is that we have
$$
\sup_{0\leq|\tau|,|\eta|\leq\theta^{-1}\leq R_0}\|D^{\alpha}a_j(\cdot, \tau, \cdot, \eta, \theta)\|_{L^{\infty}}\leq C,\ \hbox{for}\ |\alpha|\leq k_0
$$
for $j=1,2$ due to the estimates $(\ref{APP07})_{2,3}$, where $a_j$ are defined as those in Lemma 4.1 in \cite{AWXY}.
\end{pf}

\fi

\subsection{Estimates of the approximate solutions}

\iffalse

In this part we shall establish uniform estimates of the solutions
to the linearized problem arising from the Nash-Moser-H\"omander iteration scheme. Suppose that we have constructed the approximate solutions $(u^k, v^k), k=0,1,2,...,n$ of (\ref{1.1})-(\ref{1.2}), where $u^0=u^{a3}, v^0=v^{a3}$. And we build the $(n+1)-$th approximate solution $(u^{n+1}, v^{n+1})$ as follows. Set
\begin{align}
\label{LE01}
u^{n+1}=u^n+\delta u^n,\quad v^{n+1}=v^n+\delta v^n
\end{align}
where $(\delta u^n, \delta v^n)$ solves the following linear problem.
\begin{align}
\label{4.15} \left\{
\begin{array}{ll}
\mathcal{P}'_{(u_{\theta_n}^n, v_{\theta_n}^n)}(\delta u^n,
\delta v^n)=f^n,\\
\partial_x(\bar{\rho}\delta u^n)+\partial_{\eta}(\bar{\rho}\delta v^n)=0,\\
\delta u^n|_{\eta=0}=\delta v^n|_{\eta=0}=0,\quad
\lim_{\eta\rightarrow+\infty}\delta u^n=0,\\
\delta u^n|_{t=0}=0,
\end{array}
\right.
\end{align}
where
$$
u_{\theta_n}^n=S_{\theta_n}(u^{a3}+\sum_{j=0}^{n-1}\delta u^j),
\quad
v_{\theta_n}^n=S_{\theta_n}(v^{a3}+\sum_{j=0}^{n-1}\delta
v^j)
$$

\fi

To study the solutions $(\delta u^n, \delta v^n)$ to the problem \eqref{4.2} with $f^n$ given in 
\eqref{4.5}, 
as in Section 3,  set
\begin{equation}\label{omega-n}
\omega^n=\partial_{\eta}\left(\frac{\bar{\rho}\delta u^n}{\partial_{\eta}u_{\theta_n}^n}\right).
\end{equation}
Then $\omega^n$ satisfies 
\begin{align}
\label{4.16} \left\{
\begin{array}{ll}
\partial_t\omega^n+\partial_x(u_{\theta_n}^n\omega^n)+\partial_{\eta}(v_{\theta_n}^n\omega^n)-\frac{2}{\bar{\rho}}(\omega^n\chi^n)_{\eta}
+(\xi^n\int_0^{\eta}\omega^n(t,x,\tilde{\eta})d\tilde{\eta})_{\eta}-\frac{\bar{\rho}_t}{\bar{\rho}}\omega^n-\frac1{\bar{\rho}}
\omega^n_{\eta\eta}=\tilde{f}^n_{\eta},\\
\frac{1}{\bar{\rho}}(\omega^n_{\eta}+2\omega^n\chi^n)|_{\eta=0}=-\tilde{f}^n|_{\eta=0},\\
\omega^n|_{t=0}=0,
\end{array}
\right.
\end{align}
where
$$
\chi^n=\frac{\partial_{\eta}^2u^n_{\theta_n}}{\partial_{\eta}u_{\theta_n}^n},\qquad
\xi^n=\frac{(\partial_t+u^n_{\theta_n}\partial_x+v^n_{\theta_n}\partial_{\eta}-\frac1{\bar{\rho}}\partial_{\eta}^2)\partial_{\eta}u^n_{\theta_n}}
{\partial_{\eta}u_{\theta_n}^n}-\frac{u_{\theta_n}^n\bar{\rho}_x}{\bar{\rho}}\triangleq\xi^n_1-\xi^n_2,
$$
and
\begin{align}
\label{4.17}
\tilde{f}^n=\frac{\bar{\rho}f^n}{\partial_{\eta}u^n_{\theta_n}}.
\end{align}
Similar to \eqref{lambda}-\eqref{lambda-k}, we define
\begin{align*}
\lambda_{k_1,k_2}=\|u_{\theta_n}^n-u^{a3}\|_{\mathcal{B}^{k_1,k_2}_{l}}+
\|Z^{k_1}\partial_{\eta}^{k_2}v^n_{\theta_n}\|_{L^{\infty}_{\eta}(L^2_{t,x})}
+\|Z^{k_1}\partial_{\eta}^{k_2}\chi^n\|_{L^\infty_{\eta}(L^{2}_{t,x})}+\|\xi^n_1\|_{\mathcal{B}^{k_1,k_2}_{l}},
\end{align*}
and
\begin{align*}
\lambda^n_k=\sum_{k_1+[(k_2+1)/2]\leq k}\lambda^n_{k_1,k_2}.
\end{align*}

Applying Theorem \ref{Thm3.1} to the linearized problem
(\ref{4.16}), we have
\begin{thm}
\label{Thm4.1} Suppose the known functions $(\bar{\rho}, U, V)(t,x)$ satisfy the same assumptions as in Theorem \ref{MAIN}, and the main assumptions (H) are
satisfied. Then for any fixed $l>1/2$, the following estimate holds for the solution of the problem \eqref{4.16}, 
\begin{align}
\label{4.18} \|\omega^n\|_{\mathcal{A}^k_l}\leq
C_1(\lambda_4^n)\|\tilde{f}^n\|_{\mathcal{A}^k_l}+C_2(\lambda^n_4)\lambda^n_k\|\tilde{f}^n\|_{\mathcal{A}^3_l}.
\end{align}
\end{thm}
Similar to the Lemma 5.3 in \cite{AWXY}, we also have $$\|\frac{f^a}{\partial_{\eta}u^{a3}}\|_{\mathcal{A}^{k_0}_l([0,T]\times\mathbb{T}\times\mathbb{R}^+)}\leq C\varepsilon,$$
because the construction of $(u^{a3}, v^{a3})$ and the estimates in Proposition \ref{prop4.1}. Where $\varepsilon$ comes from the smallness of the integral interval of time. Then, as in \cite{AWXY}, by studying estimates of $\tilde{f}^n$ and using an induction argument, we have

\begin{thm}
\label{Thm4.2} Under the same assumptions as those in Theorem \ref{Thm4.1}, there exists a
positive constant $C_0$ such that
\begin{align}
\label{4.19} \|\omega^n\|_{\mathcal{A}^k_l}\leq
C_0\varepsilon\theta_n^{\max\{3-\tilde{k},k-\tilde{k}\}}\triangle\theta_n,
\end{align}
holds for all $n\geq 0, 0\leq k \leq k_0$ and $\tilde{k}\geq 6$ here
$\theta_n=\sqrt{\theta_0^2+n}$ and
$\triangle\theta_n=\theta_{n+1}-\theta_n$.
\end{thm}

Using the transformation \eqref{omega-n}, we can obtain

\begin{cor}
\label{Thm4.3} Under the same assumptions as those in Theorem \ref{Thm4.2}, the following estimates hold
\begin{align}
\label{PR1}
\|\delta u^n\|_{\mathcal{A}_l^k}\leq C\varepsilon\theta_j^{\max\{3-\tilde{k}, k-\tilde{k}\}}, \quad 0\leq k\leq k_0,
\end{align}
and
\begin{align}
\label{PR2}
\|\delta v^n\|_{\mathcal{D}_0^k}\leq C_1\varepsilon\theta_j^{\max\{3-\tilde{k}, k+1-\tilde{k}\}}, \quad 0\leq k\leq k_0-1.
\end{align}
\end{cor}

\subsection{Existence  to the nonlinear problem}

To show the existence of solution to the nonlinear boundary layer
equations \eqref{1.1}, we need to show the convergence of the iteration scheme \eqref{4.1}-\eqref{4.2}.  From this iteration, we know that the approximate solutions
$(u^{n+1}, v^{n+1})$ solve the following problem
\begin{align}
\label{4.32} \left\{\begin{array}{ll} \mathcal{P}(u^{n+1},
v^{n+1})=(1-S_{\theta_n})\sum_{j=0}^ne_j+S_{\theta_n}e_n+(1-S_{\theta_n})f^{a},\\
\partial_x(\bar{\rho}u^{n+1})+\partial_{\eta}(\bar{\rho} v^{n+1})=-\bar{\rho}_t,\\
u^{n+1}|_{\eta=0}=v^{n+1}|_{\eta=0}=0,\quad
\lim_{\eta\rightarrow+\infty}u^{n+1}=U(t,x),\\
u^{n+1}|_{t=0}=u_0(x,\eta).
\end{array}
 \right.
\end{align}
From the estimates given in Corollary \ref{Thm4.3}, we know
that there exist functions
$u\in\mathcal{A}_l^{\tilde{k}-2}$ and
$v\in\mathcal{D}^{\tilde{k}-3}_0$, such that $u^{n}$ converges
to $u$ in $\mathcal{A}_l^{\tilde{k}-2}$ and $v^{n}$ converges
to $v$ in $\mathcal{D}_0^{\tilde{k}-3}$. In order to show the function
pair $(u, v)$ is indeed a solution to the system (\ref{1.1}),
it suffices to show that the right hand side in equation $(\ref{4.32})_1$
converges to zero as $n$ tends to $+\infty$. Firstly, by using Lemma 4.1, 
\begin{align*}
\|(1-S_{\theta_n})(f^a+\sum_{j=0}^ne_j)\|_{\mathcal{A}_l^k}\leq
C\theta_n^{-1}(\|f^a\|_{\mathcal{A}_l^{k+1}}+\|\sum_{j=0}^ne_j\|_{\mathcal{A}_l^{k+1}}).
\end{align*}
Then it suffices to show the $\|\sum_{j=0}^ne_j\|_{\mathcal{A}_l^{k+1}}$
converges. From the definition of $e_j=e_j^1+e_j^2$ given in \eqref{e-1}-\eqref{e-2},  we have
\begin{align*}
\|e_j^1\|_{\mathcal{A}_l^{k+1}}\leq&
C(\|\delta u^j\|_{L^{\infty}}\|\delta u^j\|_{\mathcal{A}_l^{k+2}}+\|\delta
v^j\|_{L^{\infty}}\|\delta u^j\|_{\mathcal{A}_l^{k+2}}+\|\delta
v^j\|_{\mathcal{D}^{k+2}_0}\|\delta u^j\|_{L^2_{\eta,l}(L^{\infty}_{t,x})})\\
\leq& C\varepsilon^2\theta_j^{3-\tilde{k}+\max\{3-\tilde{k},
k+2-\tilde{k}\}}(\triangle\theta_j)^2
\leq C\varepsilon^2\theta_j^{k+5-2\tilde{k}}\triangle\theta_j,
\end{align*}
for $k\leq \tilde{k}-5$. And
\begin{align*}
\|e_j^2\|_{\mathcal{A}_l^{k+1}}\leq&\|(1-S_{\theta_j})(u^j-u^{a3})\partial_{\eta}(\delta
v^j)\|_{\mathcal{A}_l^{k+1}}+2\|\frac{\bar{\rho}_x}{\bar{\rho}}(1-S_{\theta_j})(u^j-u^{a3})(\delta
u^j)\|_{\mathcal{A}_l^{k+1}}\\
&+\|\partial_{\eta}((1-S_{\theta_j})(v^j-v^{a3}))(\delta
u^j)\|_{\mathcal{A}_l^{k+1}}+\|((1-S_{\theta_j})(v^j-v^{a3}))\partial_{\eta}(\delta
u^j)\|_{\mathcal{A}_l^{k+1}}\\
&+\|\partial_{\eta}((1-S_{\theta_j})(u^j-u^{a3}))(\delta
v^j)\|_{\mathcal{A}_l^{k+1}}\\
\leq& C(\|u^j-u^{a3}\|_{L^2_{\eta,l}(L^{\infty}_{t,x})}\|\delta
v^j\|_{\mathcal{D}^{k+2}_0}+\|u^j-u^{a3}\|_{\mathcal{A}_l^{k+1}}\|\partial_{\eta}(\delta
v^j)\|_{L^{\infty}}\\
&+\|u^j-u^{a3}\|_{\mathcal{A}_l^{k+1}}\|\frac{\bar{\rho}_x}{\bar{\rho}}\delta u\|_{L^{\infty}}
+\|u^j-u^{a3}\|_{L^{\infty}}\|\frac{\bar{\rho}_x}{\bar{\rho}}\delta u\|_{\mathcal{A}_l^{k+1}}\\
&+\|\delta u\|_{\mathcal{A}_l^{k+1}}\|\partial_{\eta}(v^j-v^{a3})\|_{L^{\infty}}
+\|\delta u\|_{L^2_{\eta,l}(L^{\infty}_{t,x})}\|(v^j-v^{a3})\|_{\mathcal{D}^{k+2}_0}\\
&+\|\delta u\|_{\mathcal{A}_l^{k+2}}\|(v^j-v^{a3})\|_{L^{\infty}}
+\|\partial_{\eta}\delta u\|_{L^2_{\eta,l}(L^{\infty}_{t,x})}\|(v^j-v^{a3})\|_{\mathcal{D}^{k+1}_0}\\
&+\|\delta
v^j\|_{\mathcal{D}^{k+1}_0}\|\partial{\eta}(u^j-u^{a3})\|_{L^{2}_{\eta,l}(L^{\infty}_{t,x})}
+\|\delta v^j\|_{L^{\infty}}\|u^j-u^{a3}\|_{\mathcal{A}_l^{k+2}})\\
\leq& C\varepsilon^2\theta_j^{k+3-\tilde{k}}\triangle\theta_j,
\end{align*}
for $k\leq \tilde{k}-5$. Thus, we get that
\begin{align*}
\sum_{j=0}^{+\infty}\|e_j\|_{\mathcal{A}_l^k}\leq
C\sum_{j=0}^{+\infty}\theta_j^{k+3-\tilde{k}}\triangle\theta_j\leq
CC_0,
\end{align*}
for $k\leq \tilde{k}-5$. 

Therefore, the right hand side of $(\ref{4.32})_1$
tends to zero as $n$ tends to $+\infty$. The uniqueness of
classical solutions to (\ref{1.1}) can be proved  as  in \cite{AWXY}. Then we complete the proof of Theorem \ref{MAIN}.

\vspace{.15in}

{\bf Acknowledgements:}
The first author's research was supported in part by
National Natural Science Foundation of China (NNSFC) under Grants
No. 10971134, No. 11031001 and No. 91230102. The second author is supported by NSFC
No.11171213,  Shanghai Rising Star Program No.12QA1401600. The last author's research was supported by the General Research Fund of Hong Kong,
CityU No. 103713.

\end{document}